\documentclass[reqno]{amsart}

\usepackage{amssymb,bm,mathrsfs,color}
\usepackage{hyperref}

\newtheorem{theorem}{Theorem}[section]
\newtheorem{lemma}[theorem]{Lemma}
\newtheorem{conjecture}[theorem]{Conjecture}
\theoremstyle{remark}
\newtheorem{example}[theorem]{Example}

\DeclareMathOperator{\rank}{rank}
\DeclareMathOperator{\trace}{trace}
\DeclareMathOperator{\diag}{diag}

\renewcommand{\le}{\leqslant}
\renewcommand{\ge}{\geqslant}

\renewenvironment{description}
{\begin{list}{}{%
\setlength{\labelwidth}{\leftmargin}%
\advance \labelwidth-\labelsep}}%
{\end{list}}

\begin{document}

\title{Commutative association schemes}
\author{William J. Martin}
\address{Department of Mathematical Sciences and Department of Computer Science, Worcester Polytechnic Institute, 100 Institute Road, Worcester, MA 01609, USA}
\email{martin@wpi.edu}
\author{Hajime Tanaka}
\address{Graduate School of Information Sciences, Tohoku University, Sendai 980-8579, Japan}
\email{htanaka@math.is.tohoku.ac.jp}
\date{}
\begin{abstract}
Association schemes were originally introduced by Bose and his co-workers
in the design of statistical experiments. Since that point of inception,
the concept has proved useful in the study of group actions, in algebraic
graph theory, in algebraic coding theory, and in areas as far afield as
knot theory and numerical integration. This branch of the theory, viewed
in this collection of surveys as the ``commutative case,'' has seen
significant activity in the last few decades. The goal of the present
survey is to discuss the most important new developments in several
directions, including Gelfand pairs, cometric association schemes,
Delsarte Theory, spin models and the semidefinite programming technique.
The narrative follows a thread through this list of topics, this being
the contrast between combinatorial symmetry and group-theoretic symmetry,
culminating in Schrijver's SDP bound for binary codes (based on group
actions) and its connection to the Terwilliger algebra (based on combinatorial
symmetry). We propose this new role of the Terwilliger algebra in Delsarte
Theory as a central topic for future work.
\end{abstract}
%\subjclass[2000]{05E35; 05E30}

\maketitle

\section{Introduction}

The concept of (symmetric) association schemes was first given in the design of experiments \cite{BN1939S,BS1952JASA}.
It can also be viewed as a purely combinatorial generalization of the concept of finite transitive permutation groups.\footnote{The theory of association schemes is famously said to be a ``group theory without groups'' \cite{BI1984B}.}
The Bose--Mesner algebra, which is a fundamental tool in the theory, was introduced in \cite{BM1959AMS}.
The monumental thesis of P. Delsarte \cite{Delsarte1973PRRS} proclaimed the importance of  commutative association schemes as a unifying framework for coding theory and design theory.
In \cite{Delsarte1973PRRS}, the method of linear programming was successfully combined with the duality of the (commutative) Bose--Mesner algebra, and it has been serving as one of the underlying principles in the theory of commutative association schemes.

The theory continues to grow rapidly, enlarging its diverse connections and applications to other branches of mathematics.\footnote{Bannai \cite[p.~108]{Bannai1990P} also states that ``it seems that commutative case forms a fairly closed universe, similar to the compact symmetric homogeneous spaces.''}
However, a number of important topics are either just glanced upon or not treated at all in this paper. We are limited
not only by length, but also by our own limited expertise.
Notable omitted areas include: distance-regular graphs,\footnote{\label{footnote:excuse}Indeed, a separate, comprehensive update on distance-regular graphs by different authors is reportedly in preparation and was originally planned for this collection of surveys.} the Terwilliger algebra,\footnote{While we cannot begin to cover this important topic in the present paper, we shall encounter the Terwilliger algebra, particularly in our vision for Delsarte Theory and our treatment of semidefinite programming.} and the role of association schemes in designed experiments.
In fact, if we had tried to cover all of these topics in full detail, then this survey article would have been a book rather than a paper!
However, fortunately there are a number of excellent articles/books on the theory of (commutative) association schemes and Delsarte Theory, e.g., \cite{Delsarte1973PRRS,BI1984B,BCN1989B,DL1998IEEE,Camion1998B,Bannai1999N,Bailey2004B}.
The following are a few of the  books which include accounts on commutative association schemes:
\cite{MS1977B,CL1991B,Godsil1993B,LW2001B}.
Thus, naturally guided in part by own current research interests, we shall in this paper focus  on recent progress in the theory 
that has not been treated (in detail) in previous literature. Thus our hope is to contribute an addendum to the important references
listed above; in no way do we intend to supplant any of them or summarize their content.

We now give a summary of the present paper, together with some additional comments.
\S \ref{sec:commutative-schemes} reviews the basic theory of commutative association schemes, the goal of the exposition being 
to provide  just enough background to support the later sections.
We refer the reader to \cite{BI1984B,BCN1989B} for a more comprehensive and detailed account of the theory.
In \S \ref{sec:Gelfand-pairs} we consider commutative association schemes obtained from transitive group actions.
The Bose--Mesner algebra of such an association scheme coincides with the corresponding centralizer algebra, and its eigenmatrices are equivalent to the zonal spherical functions.
We briefly introduce a beautiful phenomenon which can be observed in many families of commutative association schemes related to classical groups over finite fields, and discuss its recent application to the construction of \emph{Ramanujan graphs} which are of great interest in computer science.
The topics of \S \ref{sec:metric-cometric-schemes} are metric (or $P$-polynomial) association schemes and cometric (or $Q$-polynomial) association schemes.\footnote{The concept of metric association scheme is essentially the same as that of a distance-regular graph, which is in turn a combinatorial version of the concept of  a distance-transitive graph.}
These concepts were introduced by Delsarte \cite{Delsarte1973PRRS}, and certain systems of orthogonal polynomials naturally come into play in the theory.
The association schemes which are both metric and cometric may be viewed as finite analogues of rank one symmetric spaces \cite[pp.~311--312]{BI1984B}, and many researchers have been working on the classification of such association schemes.
Special attention will be paid, however, to the class of cometric (but not necessarily metric) association schemes, which has become an active area of research in its own right over the last decade.

\S\S \ref{sec:codes-designs}--\ref{sec:ubiquity} discuss Delsarte's theory and related topics.
We especially recommend the survey articles \cite{DL1998IEEE,Camion1998B} for the (major) progress up to 1998.
We shall see that codes and designs are dealt with in a unified manner within the framework of commutative association schemes.
While the minimum distance and dual distance (or maximum strength) are two important parameters of a code or design in Delsarte Theory, 
a similar theory exists --- with an entirely different class of examples --- for two new parameters, namely \emph{width} and \emph{dual width}, 
which were introduced  by  Brouwer et al.~in \cite{BGKM2003JCTA}.
These parameters will also be briefly reviewed in \S \ref{sec:codes-designs}.
In \S \ref{sec:duality},  we first recall standard facts on translation association schemes and their group codes.
The duality between codes and designs observed in \S \ref{sec:codes-designs} becomes more than formal in this case.
We also discuss dualities of association schemes in connection with \emph{spin model theory} (see e.g., \cite{JMN1998JAC}).
The subject of \S\ref{sec:LP-bound} is the famous linear programming bound of Delsarte. Delsarte himself explored the
specialization of this technique to metric/cometric association schemes in \cite{Delsarte1973PRRS}, and there are a number of excellent 
treatments on the impact of this  technique in coding theory (e.g., \cite{MS1977B,L1998B}).
Our update takes as its point of departure the striking 2001 result of A.\ Samorodnitsky
\cite{Samorodnitsky2001JCTA} which says that something beyond the linear programming bound will be required to resolve the most
fundamental problem in algebraic coding theory, namely the determination of the asymptotically optimal rate of a binary
block code for a communication channel with given bit-error probability.  This serves as strong motivation for what is to come in \S \ref{sec:SDP-bound}. But we also aim to extend the techniques beyond the metric/cometric cases in support of the applications
in \S \ref{sec:ubiquity}.
Delsarte's theory has been quite successful for codes/designs in metric/cometric association schemes, but the purpose of \S \ref{sec:ubiquity} is to introduce far broader applications of his theory. Working mainly on the ``eigenspace side,'' we demonstrate the great value
in extending the theory beyond the class of cometric schemes by simply listing a variety of unusual settings where 
Delsarte's theory applies. That is, we shall characterize various combinatorial objects as codes or designs in certain association schemes. Here, the linear ordering of eigenspaces
fundamental to a cometric association scheme is replaced by a partial order and certain well-known posets play a key role in the study
of designs, and in finding solutions to the linear programming bound.  For example, $(t,m,s)$-\emph{nets}  \cite{Niederreiter1987MM,Martin2007B}  --- which provide quasi-Monte Carlo
methods for numerical integration, simulation and optimization --- are closely related to the Delsarte designs in the ``ordered
Hamming scheme,'' a family of association schemes generalizing the ordinary Hamming scheme. Here the eigenspaces are
indexed by the members of a downset (or ``lower ideal'') in Young's lattice.

In \S \ref{sec:Terwilliger-algebra}, we briefly investigate the \emph{Terwilliger} (or \emph{subconstituent}) \emph{algebra} of an association scheme \cite{Terwilliger1992JAC,Terwilliger1993JACa,Terwilliger1993JACb}.
This  noncommutative matrix algebra contains much more structural information about the association scheme than the 
(commutative) Bose--Mesner algebra.
The Terwilliger algebra has proven to be a powerful tool in the study of metric \& cometric association schemes.
Despite the importance of this connection, our focus in this section is instead on the use of this algebra 
in the analysis of codes and designs, mirroring the use of  the Bose--Mesner algebra in  Delsarte Theory.
This is still quite a new approach, but we have included the account here to propose it as a possible research direction.
\S\ref{sec:SDP-bound} is devoted to a discussion of the \emph{semidefinite programming bound} introduced recently by A. Schrijver \cite{Schrijver2005IEEE}.
This new bound was first established for binary codes and constant weight codes in \cite{Schrijver2005IEEE} using the Terwilliger algebra, and then for nonbinary codes in \cite{GST2006JCTA}.
The semidefinite programming bound is always at least as good as Delsarte's linear programming bound, and numerical computations show that there are many cases where it in fact improves upon known upper bounds.
For simplicity, our exposition is restricted to binary codes.
A survey paper \cite{Vallentin2009LAA} also contains an account on this bound for binary codes based on the results in \cite{Dunkl1976IUMJ}, but we shall particularly emphasize a viewpoint centered on the Terwilliger algebra.
The current formulation of the semidefinite programming bound relies heavily on certain group actions, so that strictly speaking it does \emph{not} belong to the ``association scheme theory'' yet. However, experience shows that group actions can be supplanted with
assumptions of combinatorial regularity and 
our treatment here reflects our hope that, for a wide class of schemes, Delsarte Theory will be reconstructed entirely in the 
near future based on the Terwilliger algebra and the semidefinite programming method.

\section{Commutative association schemes}\label{sec:commutative-schemes}

Let $X$ be a finite set and $\mathbb{C}^{X\times X}$ the set of complex matrices with rows and columns indexed by $X$.
Let $\mathcal{R}=\{R_0,R_1,\dots,R_n\}$ be a set of non-empty subsets of $X\times X$.
For each $i$, let $A_i\in\mathbb{C}^{X\times X}$ be the adjacency matrix of the graph $(X,R_i)$ (directed, in general).
The pair $(X,\mathcal{R})$ is an \emph{association scheme}\footnote{The notion coincides with that of \emph{homogeneous coherent configuration}; see \cite{EP-ejc}.} \emph{with} $n$ \emph{classes} if
\begin{description}
\item[(AS1)] $A_0=I$, the identity matrix;
\item[(AS2)] $\sum_{i=0}^nA_i=J$, the all ones matrix;
\item[(AS3)] $A_i^{\mathsf{T}}\in\{A_0,A_1,\dots,A_n\}$ for $0\le i\le n$;
\item[(AS4)] $A_iA_j$ is a linear combination of $A_0,A_1,\dots,A_n$ for $0\le i,j\le n$.
\end{description}
By (AS1) and (AS4) the vector space $\bm{A}$ spanned by the $A_i$ is an algebra; this is the \emph{Bose--Mesner} (or \emph{adjacency}) \emph{algebra} of $(X,\mathcal{R})$.
The $A_i$ are linearly independent by (AS2) and thus form a basis of $\bm{A}$.
We say that $(X,\mathcal{R})$ is \emph{commutative} if $\bm{A}$ is commutative, and that $(X,\mathcal{R})$ is \emph{symmetric} if the $A_i$ are symmetric matrices.
A symmetric association scheme is commutative.
Below are listed a few examples of (symmetric) association schemes:

\begin{example}\label{exmp:Johnson}
The Johnson scheme $J(v,n)$ $(v\ge 2n)$:
$X$ is the set of all $n$-element subsets of a fixed set $\Omega$ with $v$ points, and $(x,y)\in R_i$ if $|x\cap y|=n-i$.
\end{example}

\begin{example}\label{exmp:Hamming}
The Hamming scheme $H(n,q)$ ($q\ge 2$):
$X$ is the set of all words of length $n$ over an alphabet $\mathcal{Q}$ with $q$ symbols, and $(x,y)\in R_i$ if $x$ and $y$ differ in exactly $i$ coordinate positions.
\end{example}

Let $\mathbb{C}^X$ be the set of complex column vectors with coordinates indexed by $X$, and observe that $\mathbb{C}^{X\times X}$ acts on $\mathbb{C}^X$ from the left.
For each $x\in X$ let $\hat{x}$ be the vector in $\mathbb{C}^X$ with a $1$ in coordinate $x$ and $0$ elsewhere.
We endow $\mathbb{C}^X$ with the standard Hermitian form $\langle,\rangle$ so that the $\hat{x}$ form an orthonormal basis for $\mathbb{C}^X$.

For the rest of this section, let $(X,\mathcal{R})$ be a commutative association scheme with adjacency matrices $A_0,A_1,\dots,A_n$ and Bose--Mesner algebra $\bm{A}$.
By (AS3) $\bm{A}$ is closed under conjugate transposition.
Since $\bm{A}$ is commutative, it follows that there is a unitary matrix $U\in\mathbb{C}^{X\times X}$ such that $U^{-1}\bm{A}U$ consists of diagonal matrices only.
In other words, $\mathbb{C}^X$ is decomposed as an orthogonal direct sum of $n+1$ maximal common eigenspaces\footnote{They are sometimes called the \emph{strata} \cite{Bailey2004B}.} of $\bm{A}$:
\begin{equation}
	\mathbb{C}^X=V_0\perp V_1\perp\dots\perp V_n.
\end{equation}
For each $i$, let $E_i\in\mathbb{C}^{X\times X}$ be the orthogonal projection onto $V_i$.
Then the $E_i$ form a basis of the primitive idempotents of $\bm{A}$, i.e, $E_iE_j=\delta_{ij}E_i$, $\sum_{i=0}^nE_i=I$.
Note that $|X|^{-1}J$ is an idempotent in $\bm{A}$ with rank one, hence must be primitive;
we shall always set $E_0=|X|^{-1}J$.
It also follows from (AS2) that $\bm{A}$ is closed under entrywise (\emph{Hadamard} or \emph{Schur}) multiplication, denoted $\circ$.
The $A_i$ form a basis of the primitive idempotents of $\bm{A}$ with respect to this multiplication, i.e., $A_i\circ A_j=\delta_{ij}A_i$, $\sum_{i=0}^nA_i=J$.

The \emph{intersection numbers} $p_{ij}^k$ and the \emph{Krein parameters} $q_{ij}^k$ $(0\le i,j,k\le n)$ of $(X,\mathcal{R})$ are defined by the equations
\begin{equation}
	A_iA_j=\sum_{k=0}^np_{ij}^kA_k, \quad E_i\circ E_j=\frac{1}{|X|}\sum_{k=0}^nq_{ij}^kE_k.
\end{equation}
The $p_{ij}^k$ are nonnegative integers.
On the other hand, since each $E_i\circ E_j$ (being a principal submatrix of $E_i\otimes E_j$) is positive semidefinite, it follows that the $q_{ij}^k$ are real and nonnegative.
This important restriction is known as the \emph{Krein condition}.

The change-of-basis matrices $P$ and $Q$ are defined by
\begin{equation}
	A_i=\sum_{j=0}^nP_{ji}E_j, \quad E_i=\frac{1}{|X|}\sum_{j=0}^nQ_{ji}A_j.
\end{equation}
We shall refer to $P$ and $Q$ as the \emph{first} and \emph{second eigenmatrix} of $(X,\mathcal{R})$, respectively.
Note that $P_{0i},P_{1i},\dots,P_{ni}$ give the eigenvalues of $A_i$.
The matrix $P$ is also called the \emph{character table} of $(X,\mathcal{R})$, and in fact it can be viewed as a natural generalization of the character table of a finite group;\footnote{\label{footnote:splitting-field}It is a long-standing problem whether the $P_{ij}$ are contained in a cyclotomic number field or not in general \cite[p. 123]{BI1984B}. It is known that this is the case if the $q_{ij}^k$ are rational \cite{Munemasa1991JCTA}. A negative answer to the problem would imply that the character theory of commutative association schemes is ``far'' from that of finite groups. See also \cite{Komatsu2008KJM}.} see Example \ref{exmp:group-association-scheme}.
Let
\begin{equation}
	k_i=P_{0i}, \quad m_i=Q_{0i}.
\end{equation}
It follows that $k_i$ is the valency of the regular graph $(X,R_i)$ and $m_i=\trace(E_i)=\rank(E_i)$.
The $m_i$ are called the \emph{multiplicities} of $(X,\mathcal{R})$.
For convenience, set $\Delta_k=\diag(k_0,k_1,\dots,k_n)$ and $\Delta_m=\diag(m_0,m_1,\dots,m_n)$.
Then we have
\begin{equation}
	\Delta_mP=\overline{Q}^{\mathsf{T}}\Delta_k.
\end{equation}
This is verified by evaluating $\trace(A_iE_j)$ in two ways.
Combining this with the obvious equality $PQ=QP=|X|I$, we get the following \emph{orthogonality relations}:
\begin{equation}\label{eqn:orthogonality}
	P^{\mathsf{T}}\Delta_m\overline{P}=|X|\Delta_k, \quad Q^{\mathsf{T}}\Delta_k\overline{Q}=|X|\Delta_m.
\end{equation}
We record here the eigenmatrix $P$ for Examples \ref{exmp:Johnson} and \ref{exmp:Hamming}.
See \cite{Delsarte1973PRRS,Delsarte1976JCTA,Delsarte1978SIAM,BI1984B,Stanton1985JCTA} for the details.

\begin{example}\label{exmp:dual-Hahn}
Suppose $(X,\mathcal{R})=J(v,n)$.
Then $k_i=\binom{n}{i}\binom{v-n}{i}$, $m_i=\binom{v}{i}-\binom{v}{i-1}$ and
the $P_{ij}$ are given by dual Hahn polynomials \cite[\S 1.6]{KS1998R}:
\begin{equation*}
	\frac{P_{ij}}{k_j}={}_3F_2\bigg(\!\!\!\begin{array}{c} -i,-j,i-v-1 \\ n-v, -n \end{array}\!\!\!\biggm| 1\biggr).
\end{equation*}
\end{example}

\begin{example}\label{exmp:Krawtchouk}
Suppose $(X,\mathcal{R})=H(n,q)$.
Then $k_i=m_i=\binom{n}{i}(q-1)^i$ and
the $P_{ij}$ are given by Krawtchouk polynomials \cite[\S 1.10]{KS1998R}:
\begin{equation*}
	\frac{P_{ij}}{k_j}={}_2F_1\bigg(\!\!\!\begin{array}{c} -i,-j \\ -n \end{array}\!\!\!\biggm| \frac{q}{q-1} \biggr).
\end{equation*}
\end{example}

We remark that the polynomials in Examples \ref{exmp:dual-Hahn} and \ref{exmp:Krawtchouk} belong to the so-called Askey-scheme of (basic) hypergeometric orthogonal polynomials \cite{KS1998R}, and the equations in \eqref{eqn:orthogonality} for the corresponding association schemes  amount to the orthogonality relations of these polynomials and their duals; see \S \ref{sec:metric-cometric-schemes}.

An association scheme $(X,\mathcal{S})$ on the same vertex set $X$ is called a \emph{fusion} of $(X,\mathcal{R})$ if each $S\in\mathcal{S}$ is the union of some of the $R_i$.\footnote{By (AS1), $R_0\in\mathcal{S}$. As an extreme case, we call $(X,\mathcal{R})$ \emph{amorphous} (or \emph{amorphic}) if every ``merging'' operation on $\{R_1,R_2,\dots,R_n\}$ yields a fusion; see \cite{DMpre} for a recent survey on this topic.}
As the adjacency matrices (resp. primitive idempotents) of  a fusion of $(X,\mathcal{R})$ must be $01$-linear combinations of the $A_i$ (resp. $E_i$), it is theoretically possible to find all fusions of $(X,\mathcal{R})$ from the eigenmatrix $P$.
This is accomplished using the \emph{Bannai--Muzychuk Criterion} \cite{Bannai1991JA,Muzychuk1988D}.\footnote{See e.g., \cite{Bannai1991JA,Tanaka2002EJC,Fujisaki2004DM,IMpre} for examples of explicit constructions of fusions using this criterion.}

We close the section with a brief review of subschemes and quotient schemes.
For a subset $Y\subseteq X$, define $\mathcal{R}^Y=\{R_i^Y:0\le i\le n,\ R_i^Y\ne\emptyset\}$ where we write $R^Y=R\cap(Y\times Y)$ for $R\subseteq X\times X$.
We call $(Y,\mathcal{R}^Y)$ a \emph{subscheme} of $(X,\mathcal{R})$ if it is an association scheme.
For example, $J(v,n)$ can be naturally viewed as a subscheme of $H(v,2)$.

We say that $(X,\mathcal{R})$ is \emph{primitive} if the graphs $(X,R_i)$ $(1\le i\le n)$ are connected,
and \emph{imprimitive} otherwise.
Let $I_r$ (resp. $J_r$) denote the $r\times r$ identity (resp. all ones) matrix.
Then
\begin{lemma}\label{lem:imprimitivity}
The following are equivalent:
\begin{enumerate}
\item $(X,\mathcal{R})$ is imprimitive.
\item There is a subset $\mathcal{I}\subseteq\{0,1,\dots,n\}$ such that $\sum_{i\in\mathcal{I}}A_i=I_r\otimes J_s$ for some integers $r,s\ge 2$ and an ordering of $X$.\footnote{Such a subset is often said to be \emph{closed}; see \cite{Zieschang2005B}.}
\item There is a subset $\mathcal{J}\subseteq\{0,1,\dots,n\}$ such that $\sum_{i\in\mathcal{J}}E_i=s^{-1}I_r\otimes J_s$ for some integers $r,s\ge 2$ and an ordering of $X$.
\item There exist $i\in\{1,2,\dots,n\}$ and $x,y\in X$ $(x\ne y)$ such that $E_i\hat{x}=E_i\hat{y}$.
\end{enumerate}
\end{lemma}
(To prove this, proceed e.g., in the order (i)$\Leftrightarrow$(ii)$\Leftrightarrow$(iii)$\Rightarrow$(iv)$\Rightarrow$(ii).)
Suppose now that $(X,\mathcal{R})$ is imprimitive and let a subset $\mathcal{I}$ be as in (ii) above, so that $\bigcup_{i\in\mathcal{I}}R_i$ is an equivalence relation on $X$.
Each equivalence class affords a subscheme of $(X,\mathcal{R})$.
We also have a natural structure of an association scheme on the set of all equivalence classes \cite[\S 2.9]{BI1984B}, called a \emph{quotient} (or \emph{factor}) \emph{scheme} of $(X,\mathcal{R})$.\footnote{Let $E=s^{-1}I_r\otimes J_s$ be as in (iii) above. Then the Bose--Mesner algebra of the quotient scheme is canonically isomorphic to the ``Hecke algebra'' $E\bm{A}E$ (which is also an ideal of $\bm{A}$).}
In fact, there is a concept of a ``composition factor'' of $(X,\mathcal{R})$ as well as a ``Jordan-H\"{o}lder'' theorem.
See \cite{Zieschang2005B,Zieschang-ejc} for the structure theory of (general) association schemes.
Finally, it should be mentioned that if $(X,\mathcal{R})$ is primitive and symmetric then $X$ can be viewed as a set of points on a sphere in $\mathbb{R}^{m_i}$ for each $i\in\{1,2,\dots,n\}$ in view of (iv) above. This ``spherical embedding'' is known to be quite useful; see e.g.,  
\cite[Chapter 3]{BCN1989B}, \cite{BB2006CDM}.

\section{Gelfand pairs}\label{sec:Gelfand-pairs}

Suppose that a finite group $G$ acts transitively on a finite set $X$ (from the left).
Let $\mathcal{R}=\{R_0,R_1,\dots,R_n\}$ be the set of the orbits\footnote{These orbits are also referred to as the \emph{orbitals} or the $2$-\emph{orbits} of $G$ on $X$; see \cite{EP-ejc,Weilandt1964B}.} of $G$ on $X\times X$ under componentwise action, where $R_0=\{(x,x):x\in X\}$, then $(X,\mathcal{R})$ forms an association scheme.\footnote{Association schemes of this type are sometimes called \emph{Schurian} \cite{EP-ejc}.}
Let $\pi:G\rightarrow\mathbb{C}^{X\times X}$ be the permutation representation of $G$ on $X$, i.e., $\pi(g)_{xy}=\delta_{x,gy}$ $(x,y\in X,\ g\in G)$, and observe that the Bose--Mesner algebra $\bm{A}$ of $(X,\mathcal{R})$ coincides with the centralizer (or Hecke) algebra of $\pi$:
\begin{equation}
	\bm{A}=\{M\in\mathbb{C}^{X\times X}:\pi(g)M=M\pi(g)\ \text{for all}\ g\in G\}.
\end{equation}
Hence, by Schur's lemma it follows that $(X,\mathcal{R})$ is commutative if and only if $\pi$ is multiplicity-free, i.e., if and only if $\pi$ is equivalent to a direct sum of inequivalent irreducible representations of $G$.
We note that $(X,\mathcal{R})$ is symmetric if and only if $\pi$ is multiplicity-free and each irreducible constituent is realizable in $\mathbb{R}$; this is also equivalent to the condition that the action of $G$ on $X$ is \emph{generously transitive}, i.e., for any distinct $x,y\in X$ there is an element $g\in G$ such that $gx=y$ and $gy=x$.

The $G$-set $X$ can be identified with the set of left cosets $G/K$ for some subgroup $K$ of $G$.
Note that $\pi=(1_K)^G$, where $1_K$ is the trivial representation of $K$.
For brevity, in this case we shall also let the symbol $G/K$ \emph{denote} the corresponding association scheme.
It follows that $G/K$ is primitive if and only if $K$ is maximal in $G$.
The pair $(G,K)$ is called a \emph{Gelfand pair} if $(1_K)^G$ is multiplicity-free \cite[Chapter VII]{Macdonald1995B}, i.e., if $G/K$ is commutative.
If $(G,K)$ is a Gelfand pair, then the determination of the eigenmatrices of $G/K$ is in fact equivalent to describing the zonal spherical functions of $(G,K)$ \cite[\S 2.11]{BI1984B}.

The Johnson and Hamming schemes are obtained by actions of symmetric groups and their wreath products, respectively:

\begin{example}
$J(v,n)=S_v/(S_n\times S_{v-n})$.
\end{example}

\begin{example}
$H(n,q)=(S_q\wr S_n)/(S_{q-1}\wr S_n)$.
\end{example}

Concerning finite groups, we introduce another important example:

\begin{example}[The group association scheme\footnote{E. Bannai (personal communication) pointed out that the concept of \emph{supercharacters} and \emph{superclasses} of finite groups introduced in \cite{DI2008TAMS} is equivalent to that of fusions of group association schemes (in view of the Bannai--Muzychuk Criterion).}]\label{exmp:group-association-scheme}
Suppose that $X$ is a finite group and let $G=X\times X$ act on $X$ by $(x,y)z=xzy^{-1}$.
Thus $(X,\mathcal{R})=(X\times X)/\tilde{X}$ in the above sense, where $\tilde{X}=\{(x,x):x\in X\}$.
We have $(x,y)\in R_i$ if $yx^{-1}\in C_i$, where $C_0=\{1\},C_1,\dots,C_n$ are the conjugacy classes of $X$.
In this case $\bm{A}$ is isomorphic to the center of $\mathbb{C}^X$ (viewed naturally as the group algebra of $X$) by $A_i\mapsto\sum_{x\in C_i}\hat{x}$.
Hence $(X\times X)/\tilde{X}$ is commutative and the $E_i$ are in bijection with the irreducible characters $\varepsilon_i$ of $X$.
Let $T$ be the group character table of $X$ in the usual sense, and set $\Delta_k=\diag(k_0,k_1,\dots,k_n)$, $\Delta_f=\diag(f_0,f_1,\dots,f_n)$, where $k_i=|C_i|$ and $f_i$ is the degree of $\varepsilon_i$ $(0\le i\le n)$.
Then we have\footnote{\label{footnote:idempotents}The central primitive idempotent corresponding to $\varepsilon_i$ is given by $f_i|X|^{-1}\sum_{x\in X}\overline{\varepsilon_i(x)}\hat{x}$; see e.g., \cite{CR1990B}.}
\begin{equation*}
	\Delta_fP=T\Delta_k, \quad Q=\overline{T}^{\mathsf{T}}\!\Delta_f.
\end{equation*}
Note that $(X\times X)/\tilde{X}$ is primitive if and only if $X$ is a simple group.
\end{example}

The eigenmatrices of commutative association schemes have been extensively studied in the context of spherical functions.
They are of significant interest in the theory of orthogonal polynomials (see e.g., \cite{Stanton1984P}).
See also \cite{Terras1999B,CST2008B} for a wide variety of applications.\footnote{The goal of the book \cite{Terras1999B} is to consider finite analogues of the symmetric spaces including $\mathbb{R}^n$ and the Poincar\'{e} upper half plane, partly in order to ``develop an understanding of the continuous theory by developing the finite model'' \cite[p.~1]{Terras1999B}.}
Besides, there are many other commutative association schemes whose eigenmatrices possess quite beautiful structure.\footnote{Concerning the character theory of finite groups, we especially refer to \cite{BKS1990JA,Henderson2001RT} for a strong analogy between the character tables of $GL(n,q)$, $GU(n,q^2)$ and the eigenmatrices of $GL(2n,q)/Sp(2n,q)$, $GL(n,q^2)/GL(n,q)$, $GL(n,q^2)/GU(n,q^2)$. See also \cite{Bannai1990P}.}
As a typical example we consider the association schemes $O_{2m+1}(q)/O_{2m}^{\pm}(q)$, where for brevity we assume that $q$ is even.\footnote{These association schemes arise from the action of $O_{2m+1}(q)$ on each of the sets of plus-type and minus-type hyperplanes. See \cite{BHS1990JCTA} for the preceding results in the case of odd $q$.}
We first recall that $O_3(q)/O_2^-(q)$ is symmetric with $q/2-1$ classes (cf. \cite{Tanaka2001M,HX2006JAC}).
The first eigenmatrix $P$ is of the following form:\footnote{It is also known that all nontrivial multiplicities coincide (and equal $q+1$), i.e., $O_3(q)/O_2^-(q)$ is \emph{pseudocyclic} \cite[\S 2.2B]{BCN1989B}. Pseudocyclic association schemes can be used to construct strongly regular graphs as well as distance-regular graphs with diameter three; see \cite[\S 12.7]{BCN1989B}.}
\begin{equation*}
	P=\begin{bmatrix} 1 & q+1 & \cdots & q+1 \\ 1 \\ \vdots && P_0 \\ 1 \end{bmatrix},
\end{equation*}
where $P_0$ is a square matrix of size $q/2-1$.
(Recall that the top row of $P$ gives the valencies.)
Next, let $\epsilon\in\{+,-\}$ and suppose $m>1$ if $\epsilon=-$.
Then $O_{2m+1}(q)/O_{2m}^{\epsilon}(q)$ is symmetric with $q/2$ classes and the nontrivial valencies are given by $k_1=(q^{m-1}+(\epsilon 1))(q^m-(\epsilon 1))$, $k_2=k_3=\dots=k_{q/2}=q^{m-1}(q^m-(\epsilon 1))$.
Moreover it turns out that the first eigenmatrix of $O_{2m+1}(q)/O_{2m}^{\epsilon}(q)$ contains $-(\epsilon q^{m-1})P_0$ as its lower-right submatrix.\footnote{The proof is based on comparing the intersection numbers of these association schemes. Note that the other entries are determined from the orthogonality relations \eqref{eqn:orthogonality}.}
In this situation, we say that the eigenmatrix of $O_{2m+1}(q)/O_{2m}^{\epsilon}(q)$ is \emph{controlled} by that of the smaller association scheme $O_3(q)/O_2^-(q)$.
There are a lot of fascinating examples of this kind; see \cite{Bannai1990P,ST2006P} and the references therein.

We close this section with an application of this theory.
A connected $k$-regular graph is called \emph{Ramanujan} if all eigenvalues $\theta$ such that $|\theta|\ne k$ satisfy $|\theta|\le 2\sqrt{k-1}$.
Ramanujan graphs are good expanders and there are broad applications in computer science.
Moreover, these graphs are precisely the regular graphs whose Ihara zeta functions satisfy the Riemann hypothesis.
See \cite{Terras1999B} for the details.

Observe that $O_3(q)/O_2^-(q)=PGL(2,q)/D_{q+1}$ where $D_{q+1}$ is a dihedral subgroup of order $2(q+1)$, so that it is a quotient scheme of $GL(2,q)/GL(1,q^2)$.
The latter association scheme (for both even and odd $q$) is symmetric with $q-1$ classes, and its $q-2$ graphs with valency $q+1$ are called the \emph{finite upper half plane graphs} \cite{Terras1999B}.
These graphs are shown to be Ramanujan,\footnote{The proof amounts to estimating several exponential sums over finite fields, such as \emph{Soto-Andrade sums}, based on the work of A. Weil, N. Katz, W. Li and many others. See \cite{Terras1999B}.} from which it follows that the $(q+1)$-regular graphs attached to $O_3(q)/O_2^-(q)$ are Ramanujan as well.
An implication of the above comments on the eigenmatrices is that \emph{all} graphs with valency $q^{m-1}(q^m+1)$ associated with $O_{2m+1}(q)/O_{2m}^-(q)$ are Ramanujan;\footnote{Using a simple number-theoretic argument, it is also shown that for each fixed $q$ there are infinitely many values of $m$ such that the graphs attached to $O_{2m+1}(q)/O_{2m}^+(q)$ are Ramanujan.} see \cite{BST2004EJC}.

This systematic construction of Ramanujan graphs is an interesting application of the approach from the eigenmatrices of commutative association schemes.\footnote{Note, however, that this construction yields only finitely many Ramanujan graphs for each fixed valency.}
In fact, the same method works for many other examples of controlling association schemes; see \cite{BSTpre,BST2004EJC}.
See also \cite{LM2005FFA,DLM2007DM} for related constructions of Ramanujan graphs and \cite{Vinh2008EJC} for an application of the results in \cite{BSTpre,BST2004EJC} to the Erd\H{o}s distance problem.

\section{Metric/cometric association schemes}\label{sec:metric-cometric-schemes}

Suppose that $(X,\mathcal{R})$ is a symmetric association scheme.
We say that $(X,\mathcal{R})$ is \emph{metric} (or $P$-\emph{polynomial}) with respect to the ordering $\{A_i\}_{i=0}^n$ if for each $i$ $(0\le i\le n)$ there is a polynomial $v_i$ with degree $i$ such that $P_{ji}=v_i(P_{j1})$ $(0\le j\le n)$.
Such an ordering is called a $P$-\emph{polynomial ordering}.
Dually, we say that $(X,\mathcal{R})$ is \emph{cometric} (or $Q$-\emph{polynomial}) with respect to the ordering $\{E_i\}_{i=0}^n$ if for each $i$ $(0\le i\le n)$ there is a polynomial $v_i^*$ with degree $i$ such that $Q_{ji}=v_i^*(Q_{j1})$ $(0\le j\le n)$.
Such an ordering is called a $Q$-\emph{polynomial ordering}.
Note that in each of the above definitions the $v_i$ (resp. $v_i^*$) form a system of orthogonal polynomials by \eqref{eqn:orthogonality}.
Note also that $(X,\mathcal{R})$ is metric (resp. cometric) with respect to the above ordering if and only if for all $i,j,k$ $(0\le i,j,k\le n)$ we have $p_{ij}^k=0$ (resp. $q_{ij}^k=0$) if $i+j>k$ and $p_{ij}^k\ne 0$ (resp. $q_{ij}^k\ne 0$) if $i+j=k$.

A connected undirected graph $(X,R)$ with diameter $n$ and path-length distance $\partial$ is called \emph{distance-regular} if the $n+1$ relations $R_i=\{(x,y)\in X\times X:\partial(x,y)=i\}$ $(0\le i\le n)$ define an association scheme \cite{BI1984B,BCN1989B}.
Thus metric association schemes, with specified $P$-polynomial ordering, are in bijection with distance-regular graphs.
We refer the reader to \cite{BCN1989B} for the basic theory on this topic,\footnote{See also ``Additions and corrections'' to the book \cite{BCN1989B} available at A. E. Brouwer's webpage: \href{http://www.win.tue.nl/~aeb/}{http://www.win.tue.nl/\~{}aeb/} (cf. footnote \ref{footnote:excuse}).} and briefly comment on metric association schemes which are also cometric.
(Henceforth we will use the phrase ``metric \& cometric'' to describe such association schemes.)
This class contains $J(v,n)$, $H(n,q)$ and many other important examples which arise from finite classical groups and classical forms over finite fields, e.g., Grassmann schemes and bilinear forms schemes;\footnote{These are $q$-analogues of $J(v,n)$ and $H(n,q)$, respectively.} see \cite[\S 3.6]{BI1984B}, \cite[Chapter 9]{BCN1989B} and \cite{DK2005IM}.
The famous theorem of Leonard \cite{Leonard1982SIAM} states that in this case the above polynomials $v_i$ and $v_i^*$ belong to the terminating branch of the Askey scheme \cite{KS1998R} (up to normalization),  so that they are $q$-Racah polynomials \cite[\S 3.2]{KS1998R} or their special/limiting cases;\footnote{Note that we also allow the specialization $q\rightarrow -1$. Leonard's theorem was obtained just shortly after the 1979 discovery of the $q$-Racah (or Askey--Wilson) polynomials, and led Andrews and Askey \cite{AA1985P} to their definition of the \emph{classical} orthogonal polynomials. Bannai \cite[p.~27]{Bannai1999N} states that ``it is very interesting that combinatorics gave a meaningful influence to the theory of orthogonal polynomials in this way.''} see also \cite[\S 3.5]{BI1984B}.
Recently, Leonard's theorem has been reformulated in the purely linear algebraic framework of \emph{Leonard pairs} \cite{Terwilliger2001LAA,Terwilliger2004LAA}.
Leonard pairs are used to describe certain irreducible modules for the \emph{Terwilliger algebra} (\S \ref{sec:Terwilliger-algebra}) of metric \& cometric association schemes.\footnote{Some problems on metric \& cometric association schemes can be unified and most elegantly stated in terms of Leonard pairs; see \cite{Tanakapre}.}
We also remark that Leonard pairs arise naturally in other various contexts, such as in representation theory;\footnote{For example, we may obtain Leonard pairs from the finite dimensional irreducible modules for the Lie algebra $\mathfrak{sl}_2$ as well as the quantum algebra $U_q(\mathfrak{sl}_2)$.} see \cite{Terwilliger2003JCAM,Terwilliger2006N} for the details.

Compared with metric association schemes, the systematic study of cometric (but not necessarily metric) association schemes has begun rather recently.
As we shall discuss below, they are of particular interest because of their connections, e.g., to spherical designs, (Euclidean) lattices and also mutually unbiased bases in quantum information theory.
See also \cite{BB-ejc}.

The famous \emph{Bannai--Ito Conjecture} \cite[p.~237]{BI1984B} states that there are only finitely many distance-regular graphs with any given valency $k>2$ (the polygons are all distance-regular with $k=2$).
For recent activity on this conjecture, see e.g., \cite{KM2004JAC,BKM2007EJC} and the references therein.\footnote{As of this writing,  it has been announced that the full conjecture has been proven by Bang, Koolen and Moulton.}
The dual to this conjecture is the following:

\begin{theorem}[\cite{MW2009EJC}]
For each fixed $m>2$, there are only finitely many cometric association schemes with $Q$-polynomial
ordering $\{E_i\}_{i=0}^n$ satisfying $\rank(E_1)=m$.
\end{theorem}

A key step in the proof of this theorem is to bound the degree of the splitting field, based on the results of \cite{Suzuki1998JACb}; see also \cite[\S 3]{CS2009EJC}.
The \emph{splitting field} of $(X,\mathcal{R})$ is the smallest extension of the rational number field $\mathbb{Q}$ which contains all the $P_{ij}$.\footnote{See footnote \ref{footnote:splitting-field}.}
While most distance-regular graphs with classical parameters have rational splitting field, the regular $n$-gon has splitting field $\mathbb{Q}(\zeta)$ where $\zeta=\cos (2\pi/n)$ which, when $n$ is prime for example, gives a degree $(n-1)/2$ extension of $\mathbb{Q}$.
In the case of distance-regular graphs of valency $k>2$, only one known example --- the Biggs--Smith graph --- has an eigenvalue not belonging to a quadratic extension of $\mathbb{Q}$.  To push this a bit further, to our knowledge, 
the only distance-regular graphs known with splitting field not contained in a quadratic extension of $\mathbb{Q}$ are 
\begin{itemize}
\item the Biggs--Smith graph, with spectrum
\begin{gather*}
	\textstyle
	3^1, \ 2^{18}, \ 0^{17},  \  \left( \frac{1+\sqrt{17}}{2} \right)^{9}, \  
	\left( \frac{1-\sqrt{17}}{2} \right)^{9}, \
	\left(-1+ 2\cos \frac{\pi}{9}  \right)^{16},  \\
	\textstyle
	\left( -1-\cos \frac{\pi}{9}  + \sqrt{3}  \sin \frac{\pi}{9}  \right)^{16}, \ 
	\left( -1-\cos \frac{\pi}{9}  - \sqrt{3}  \sin \frac{\pi}{9}  \right)^{16}
\end{gather*}

\item incidence graphs of generalized hexagons $GH(q,q)$, with eigenvalues
$$ \pm (q+1), \qquad  0, \qquad \pm \sqrt{q}, \qquad \pm \sqrt{3q}.$$
(These can be viewed as point graphs of thin generalized $12$-gons of order $(1,q)$.)
\item line graphs of the above graphs, with eigenvalues
$$ 2q, \qquad q-1, \qquad -2,  \qquad q-1 \pm \sqrt{q},   \qquad q-1 \pm \sqrt{3q} ,
$$
which are the  point graphs of generalized $12$-gons of order $(q,1)$.
\end{itemize}

The question arises as to whether there exists a bound on the dimension of the splitting field of a distance-regular graph in terms of its valency.
Any such result would immediately imply the Bannai--Ito Conjecture since the interval $[-k,k]$ would then contain a limited number of potential eigenvalues (since it must contain all their conjugates).
But such a bound seems beyond our reach at this point.
In \cite[p.~130]{BCN1989B}, Brouwer, et al.~ask if a distance-regular graph with $k>2$ must have at least two integral eigenvalues.
Even this apparently simple question remains unresolved to date.

\begin{conjecture}[Bannai and Ito {\cite[p.~312]{BI1984B}}]
For $n$ sufficiently large, a primitive association scheme with $n$ classes is metric if and only if it is cometric.
\end{conjecture}

While no counterexamples are known to this conjecture, there are imprimitive distance-regular graphs --- the doubled Odd graphs --- which are not cometric yet have arbitrarily large diameter, and there are some $Q$-bipartite doubles of certain Hermitian forms dual polar spaces \cite[p.~315]{BI1984B} which are cometric yet not metric.
So the ``primitive'' condition in the conjecture is necessary.
Still, no progress has been made toward proving the conjecture. From the ``cometric viewpoint,'' 
perhaps the following questions will prove easier to attack:

\medskip\noindent
\textit{Question:}
Suppose that $(X,\mathcal{R})$ is cometric with $Q$-polynomial ordering $\{E_i\}_{i=0}^n$ and the $A_i$ ordered so that $Q_{01}>Q_{11} > \cdots > Q_{n1}$.
If $n$ is sufficiently large,  must the adjacency matrix $A_1$ of relation $R_1$ have $n+1$ distinct eigenvalues?

\medskip\noindent
\textit{Question:}
Does there exist an absolute constant $\nu$ such that, for any cometric association scheme $(X,\mathcal{R})$ with  $Q$-polynomial ordering $\{E_i\}_{i=0}^n$ and the $A_i$ ordered so that $Q_{01}>Q_{11} > \cdots > Q_{n1}$, we have $p_{1j}^k = 0$ whenever $k > j+\nu$?

\medskip
It is easy to prove that the valencies of any distance-regular graph with diameter $n$ satisfy the unimodal property:
\begin{equation}
	k_0  = 1 < k_1 \le k_2 \le \dots \le k_{n'} \ge k_{n'+1} \ge \dots \ge k_n
\end{equation}
for some $1\le n'\le n$, possibly $n'=n$.
For cometric association schemes, we have

\begin{conjecture}[Bannai and Ito {\cite[p.~205]{BI1984B}}]
If $(X,\mathcal{R})$ is cometric with $Q$-polynomial ordering $\{E_i\}_{i=0}^n$, then the $m_i$ satisfy the unimodal property:
\begin{equation*}
	m_0  = 1 < m_1 \le m_2 \le \dots \le m_{n'} \ge m_{n'+1} \ge \dots \ge m_n
\end{equation*}
for some $1\le n' \le n$.
\end{conjecture}

In personal communication with P. Terwilliger,  the following stronger claim was made for all cometric association schemes:

\begin{conjecture}[D. Stanton]
For $i< n/2$,  $m_i \le m_{i+1}$ and $m_i \le m_{n-i}$.
\end{conjecture}

Stanton's conjecture has been proven under the added assumption that the association scheme is \emph{dual thin} 
(\S \ref{sec:Terwilliger-algebra}) \cite{SC2001EJC} or metric \cite{Pascasio2002EJC}.

It is well known that a metric association scheme admits at most two $P$-polyn\-mi\-al 
orderings \cite[Theorem~4.2.12]{BCN1989B}.  \S 4.2D in \cite{BCN1989B} examines the 
possbilities for a second $P$-polynomial ordering of a distance-regular graph and 
obtains substantial parameter conditions.
In the cometric case, we have the following result of Suzuki (1998):

\begin{theorem}[\cite{Suzuki1998JACb}]
If $\{E_i\}_{i=0}^n$ is a $Q$-polynomial ordering for a cometric
association scheme $(X,\mathcal{R})$, then any second 
such ordering must be one of: 
\begin{itemize}
\item $E_0 , E_2 , E_4 ,  \ldots , E_3 , E_1 $;
\item $E_0 , E_n , E_1 , E_{n-1} , \ldots  $;
\item  $E_0 , E_n , E_2 , E_{n-2 }, \ldots $ ($n$ odd);
\item   $E_0, E_{n-1} , E_2 , E_{n-3} ,  \ldots$ ($n$ even);
\item $E_0 , E_5 , E_3 , E_2 , E_4 , E_1 $.
\end{itemize}
\end{theorem}

Further conditions were given in \cite{Suzuki1998JACb} (e.g., $(X,\mathcal{R})$ must be almost $Q$-bipartite,\footnote{A cometric association scheme
with $n$ classes is {\em almost $Q$-bipartite} if its parameters satisfy $a_i^*=0$ for all $i<n$, yet $a_n^* > 0$.} in the first case).
It is possible that the last case may be ruled out.

The Krein parameters of a cometric scheme $(X,\mathcal{R})$ are entirely determined by its \emph{Krein array}
\begin{equation}
	\iota^*(X,\mathcal{R})=\{b_0^*,b_1^*,\dots,b_{n-1}^*;c_1^*,c_2^*,\dots,c_n^*\},
\end{equation}
where $b_i^*=q_{1,i+1}^i$ $(0\le i\le n-1)$ and $c_i^*=q_{1,i-1}^i$ $(1\le i\le n)$.
We also define $a_i^*=q_{1i}^i$ $(0\le i\le n)$.
It is well known \cite[p.~315]{BI1984B}  that an imprimitive distance-regular graph with valency $k>2$ is bipartite or antipodal (or both).
The dual situation is not yet fully resolved.
In 1998, Suzuki proved

\begin{theorem}[\cite{Suzuki1998JACa}]
If $(X,\mathcal{R})$ is an imprimitive cometric association scheme with $Q$-polynomial ordering $\{E_i\}_{i=0}^n$, then at least one of the following holds:\footnote{With the notation of Lemma \ref{lem:imprimitivity}, the types (i)--(iv) correspond to $\mathcal{J}=\{0,2,4,\dots\}$, $\mathcal{J}=\{0,n\}$, $\mathcal{J}=\{0,3\}$ and $\mathcal{J}=\{0,3,6\}$, respectively.}
\begin{enumerate}
\item $(X,\mathcal{R})$ is $Q$-bipartite: $a_i^*=0$ for $1\le i\le n$;
\item $(X,\mathcal{R})$ is $Q$-antipodal: $b_i^*=c_{n-i}^*$ for $1\le i\le n$, except possibly $i=\lfloor n/2\rfloor$;
\item $n=4$ and $\iota^*(X,\mathcal{R})=\{m,m-1,1,b_3^*;1,c_2^*,m-b_3^*,1\}$, where $a_2^*>0$;
\item $n=6$ and $\iota^*(X,\mathcal{R})=\{ m,m-1,1,b_3^*, b_4^*, 1; \ 1, c_2^*, m-b_3^*, 1, c_5^*, m\}$, where $a_2^*= a_4^* + a_5^*> 0$.
\end{enumerate}
\end{theorem}

Schemes of  type (iii) in the theorem have  recently been ruled out \cite{CS2009EJC}.
No examples are known of type (iv); it is quite possible that none exist and then the theorem exactly mirrors the
result for imprimitive distance-regular graphs.

Let us briefly review the known examples of such imprimitive ``polynomial schemes'' with three or four classes.
In the metric case, a bipartite distance-regular graph of diameter three is necessarily the incidence graph of some symmetric $(v,k,\lambda)$ block design.
Any such scheme is cometric as well.
An antipodal distance-regular graph of diameter three is a cover of a complete graph \cite{GH1992JCTB}; these are cometric precisely when the cover has index two.
The only distance-regular graphs of diameter three which are both bipartite and antipodal are the complete bipartite graphs with a perfect matching deleted.
These trivial examples are both metric and cometric.

The 3-class imprimitive cometric schemes follow a landscape dual to this.
The $Q$-bipartite examples are all Taylor graphs; they are all index two distance-regular antipodal covers of the complete graphs.
A 3-class $Q$-antipodal scheme is equivalent to a linked system of symmetric designs \cite{CS1973IM,MMW2007JAC}; these are only metric when there are two $Q$-antipodal classes, these being the incidence graphs of symmetric designs mentioned above.
The only examples which are both $Q$-bipartite and $Q$-antipodal are again the complete bipartite graphs with a one-factor deleted.

In the case of  imprimitive $4$-class schemes, the bipartite distance-regular graphs of diameter four are incidence graphs of various designs and geometries (e.g., generalized quadrangles) while the antipodal
distance-regular graphs of diameter four are antipodal covers of strongly regular graphs.
The distance-regular graphs of diameter four which are both bipartite and antipodal are characterized as incidence graphs of symmetric $(m,\mu)$-nets \cite[p.~18]{BCN1989B}.
On the cometric side the $4$-class schemes which are $Q$-bipartite correspond to systems of lines with two angles, one of which is $\pi/2$; the $4$-class $Q$-antipodal schemes are roughly the linked systems of strongly regular designs.
Interestingly, the $4$-class schemes which are both $Q$-bipartite and $Q$-antipodal are in one-to-one correspondence with sets of \emph{real mutually unbiased bases} \cite{LMOpre}, which we now define. (See also \cite{ABS2009JCTA}.)

Let $\{b_1,b_2,\dots,b_d\}$ and $\{b'_1,b'_2,\dots,b'_d\}$ be two orthonormal bases for $\mathbb{C}^d$. We say these bases
are \emph{unbiased} (relative to one another) if $|\langle b_i | b'_j\rangle | = 1/\sqrt{d}$
for all $1\le i,j\le d$ where $\langle | \rangle$ is the standard Hermitian inner product on $\mathbb{C}^d$. A collection of orthonormal bases
for $\mathbb{C}^d$ is \emph{mutually unbiased} if any two distinct bases from the set are unbiased relative to one another. 
For $d$ a prime power, there is a construction of $d+1$ mutually unbiased bases (MUBs) in $\mathbb{C}^d$. For other dimensions
this is mostly an open question. Such constructions are useful for several applications in quantum information theory,
such as quantum key distribution and quantum state tomography. (See \cite{BSTWpre} and the references therein.)

When we restrict the bases to be real, the problem changes qualitatively. Indeed, for unit vectors $b$ and $b'$ from distinct
bases, we must have $\langle b| b'\rangle = \pm 1/\sqrt{d}$. Let $M_d$ denote the maximum possible number of real MUBs in
dimension $d$. It follows from an old result of Delsarte et al.~\cite{DGS1975PRR} that $M_d \le d/2+1$; see also \cite{CCKS1997PLMS}. This bound is achieved
for $d=4^k$ via a construction using Kerdock sets. In \cite{BSTWpre}, it is established for example that 
\begin{itemize}
\item $M_d=1$ unless $d=2$ or $4 | d$;
\item $M_d \ge 2$ if an only if there exists a Hadamard matrix of side $d$;
\item $M_d \le 3$ unless $d/4$ is an even square.
\end{itemize}
Using the results of \cite{LMOpre}, each of these results gives either a construction or a non-existence result for cometric association schemes with four classes which are both $Q$-bipartite and $Q$-antipodal; specifically, $M_d$ is an upper bound on the 
number $k$ of $Q$-antipodal classes in such association scheme on $2kd$ vertices with $Q$-antipodal classes of size $2d$.

If $(X,\mathcal{R})$ is $Q$-bipartite with $Q$-polynomial ordering 
$\{E_i\}_{i=0}^n$, then the set $\{E_1\hat{x}:x\in X\}$ is closed under multiplication by $-1$; so, viewed as points on the unit sphere in $\mathbb{R}^{m_1}$, these schemes are really best viewed as sets of lines through the origin. The imprimitivity system here has all
equivalence classes of size two and the quotient scheme, on $|X|/2$ vertices, is often interesting. Examples include the
schemes arising from the shortest vectors in the $E_6$, $E_7$, $E_8$ and Leech lattices, as well as an overlattice of 
the Barnes--Wall lattice in $\mathbb{R}^{16}$; these have Krein arrays
\begin{itemize}
\item $\iota^*(E_6) = \left\{ 6, 5,  \frac{ 9}{2 }, 3; \,  1,  \frac{ 3}{ 2}, 3, 6
\right\}$
\item $\iota^*(E_7) = \left\{ 7, 6,  \frac{ 49}{ 9},  \frac{ 35}{ 8};   \, 1,  \frac{ 14}{ 9},  \frac{ 21}{ 8}, 7
\right\}$
\item $\iota^*(E_8) = \left\{ 8,7, \frac{ 32}{ 5},6;   1, \frac{ 8}{ 5},2,8
\right\}$
\item $\iota^*(\mathit{Leech}) = \left\{ 24, 23,  \frac{ 288}{ 13},  \frac{ 150}{ 7},  \frac{ 104}{ 5},  \frac{ 81}{ 4};   \,
1,  \frac{ 24}{ 13},  \frac{ 18}{ 7},  \frac{ 16}{ 5},  \frac{ 15}{ 4}, 24
\right\}$
\item $\iota^*(\mathrm{OBW16}) = \left\{ 16, 15,  \frac{ 128}{ 9}, 8;   1,  \frac{ 16}{ 9}, 8, 16 
\right\}$
\end{itemize}
Further $Q$-bipartite examples come from $Q$-bipartite doubles of certain strongly regular graphs, such as the
two subconstituents of the McLaughlin graph.

Concerning the structure of $Q$-antipodal schemes, again very little is known. The quotient scheme is a one-class scheme.
With the natural ordering on the $A_i$,
we have $\mathcal{I}=\{0,2,4,\dots\}$ in the notation of Lemma \ref{lem:imprimitivity}.
The following theorem
has been referred to as the ``Dismantlability Theorem'':\footnote{This theorem is formally dual to an unpublished result of C. Godsil (personal communication) which states that, in an antipodal distance-regular graph, any subset of an antipodal class is a \emph{completely regular code} (\S \ref{sec:codes-designs}).}

\begin{theorem}[\cite{MMW2007JAC}]
If $(X,\mathcal{R})$ is $Q$-antipodal and $Y \subseteq X$ is a union of $\ell$ $Q$-antipodal classes,
then $(Y,\mathcal{R}^Y)$
is a cometric subscheme, which is $Q$-antipodal as well,
provided $\ell>1$.
\end{theorem}

In \S\ref{sec:duality},  we shall investigate duality among association schemes. As a special case, if $(X,\mathcal{R})$
is the coset scheme of an additive \emph{completely regular} (\S \ref{sec:codes-designs}) code $Y$ in $H(n,q)$, then the dual of $(X,\mathcal{R})$ (induced on the dual code $Y^{\circ}$) is a cometric subscheme inside $H(n,q)$.
In this way, we obtain a number
of cometric schemes from the perfect binary and ternary Golay codes and some codes derived from them \cite[p.~356]{BCN1989B}.
What is new here is that, since several of these coset graphs are antipodal, their dual schemes are $Q$-antipodal and the above
theorem gives us new cometric schemes which are not metric. We give two examples here.

\begin{itemize}
\item The dual scheme of the coset graph of the shortened ternary Golay code is a $Q$-antipodal scheme on 243 vertices
with three $Q$-antipodal classes.
If we dismantle this, taking two of these classes only, we obtain a $Q$-antipodal scheme with Krein array
$\{ 20, 18, 3, 1;   1, 3, 18, 20 \}$. It is interesting to note that the dual parameter set remains open for a possible antipodal
diameter four distance-regular graph.
\item Example (A16) on p.~365 of \cite{BCN1989B} is the coset graph of  an additive binary code derived from the extended
binary Golay code. Its dual scheme has $2048$ vertices and four $Q$-antipodal classes. If we take only three of 
these, we obtain a $Q$-antipodal scheme with Krein array
$\{ 21, 20, 16, 8, 2, 1;   1, 2, 4,$ $16, 20, 21 \}$. In this case, the dual parameter set has been shown to  be unrealizable as a distance-regular graph by counting hexagons in such a graph \cite[p.~365]{BCN1989B}.
\end{itemize}

In a terse summary of spherical designs \cite{Munemasa2007B}, Munemasa gives numerous examples of cometric schemes arising from lattices which are not distance-regular graphs.
Martin et al.~\cite{MMW2007JAC} build on this list, including some schemes coming from error-correcting codes, block designs and the above theorem applied to known $Q$-antipodal schemes.
Higman's paper \cite{Higman1995EJC}  on strongly regular designs contains further examples.

We have already mentioned some imprimitive examples. It is remarkable that very few primitive cometric association schemes are
known which are not metric. The only known examples, to our knowledge, are the following:\footnote{See also an on-line table of cometric association schemes which are not metric, maintained by W. J. Martin: \href{http://users.wpi.edu/~martin/RESEARCH/QPOL/}{http://users.wpi.edu/$\sim$martin/RESEARCH/QPOL/}}
\begin{itemize}
\item  the block scheme of the $4$-$(11,5,1)$ Witt design, with $n=3$, $|X|=66$ and Krein array
$\{ 10, 242/27, 11/5;  \ 1, 55/27, 44/5 \}$
\item the block scheme of the $5$-$(24,8,1)$ Witt design, with $n=3$, $|X|=729$ and Krein array
$\{ 23, 945/44, 1587/80;   \ 1, 345/176, 207/20 \}$
\item a spherical design derived from the Leech lattice with $n=3$, $|X|=2025$ and Krein array
$\{22,21,625/33;   \ 1, 11/6, 30/11\}$
\item the block scheme of a 4-(47,11,48) design arising from codewords of weight 11 in a certain quadratic residue
code of length $47$, with $n=3$, $|X|=4324$ and Krein array
$\{ 46, 77315/1782 , 24863/847; \   1, 37835/19602 , 2162/231 \}$
\item  the ``antipodal'' quotient of the association scheme on shortest vectors of the Leech lattice, with
$n=3$, $|X|=98280$ and Krein array \\
$\{299, 1800/7, 4563/20;$  $1, 156/35, 195/4 \}$
\item a spherical design derived from the Leech lattice
with $n=4$, $|X|=7128$ and Krein array
$\{ 22, 21, 121/6, 2187/125;  \ 1, 11/6, 363/125, 6 \}$
\item  another derived spherical design arising among the shortest vectors
of the Leech lattice, with $n=5$, $|X|=47104$ and Krein array\\
$\{ 23, 22, 529/25, 184/9, 483/25;$  $1, 46/25, 23/9, 92/25, 23/3 \}$
\end{itemize}

\section{Codes and designs}\label{sec:codes-designs}

Suppose that $(X,\mathcal{R})$ is a commutative association scheme.
Throughout this section, let $Y$ be a nonempty subset of $X$ with $1<|Y|<|X|$.
Let $\chi=\sum_{x\in Y}\hat{x}$ be the characteristic vector of $Y$.
The \emph{inner distribution} of $Y$ is the vector $\bm{a}=(a_0,a_1,\dots,a_n)$ defined by
\begin{equation}
	a_i=\frac{1}{|Y|}\chi^{\mathsf{T}}A_i\chi=\frac{1}{|Y|}|R_i\cap(Y\times Y)|.
\end{equation}
Note that the $a_i$ are nonnegative, $a_0=1$ and $(\bm{a}Q)_0=|Y|$.\footnote{In general, for a vector $\bm{c}=(c_0,c_1,\dots,c_n)$ we call $\bm{c}Q$ the \emph{MacWilliams transform} of $\bm{c}$.}
Since $(\bm{a}Q)_i=|X||Y|^{-1}\chi^{\mathsf{T}}E_i\chi$ it follows that the $(\bm{a}Q)_i$ are also \emph{real} and \emph{nonnegative}; this 
simple fact underlies Delsarte's linear programming method; see \S \ref{sec:LP-bound}.
We remark that $(\bm{a}Q)_i=0$ if and only if $E_i\chi=0$.

For a subset $\mathcal{T}$ of $\{1,2,\dots,n\}$, we call $Y$ a $\mathcal{T}$-\emph{code} (resp. (\emph{Delsarte}) $\mathcal{T}$-\emph{design}) if $a_i=0$ (resp. $(\bm{a}Q)_i=0$) for all $i\in\mathcal{T}$.
A $\{1,2,\dots,t\}$-design is simply called a $t$-\emph{design}.\footnote{\label{footnote:ordering} In what follows, if we define a concept/parameter which depends on the ordering of the $A_i$ or the $E_i$ (such as a $t$-design) then we shall understand that such an ordering is implicitly fixed. Whenever we state a result involving these concepts/parameters, the orderings will be explicitly specified or clear from the context.}
Codes in $H(n,q)$ are the familiar ``block codes of length $n$,'' and codes in $J(v,n)$ are precisely the  binary constant-weight codes.
We remark that codes in the bilinear forms schemes also have applications to \emph{space-time codes}; see \cite{GV2008AMC}.

If $Y$ is a $\mathcal{T}$-code and if $Z\subseteq X$ is a $\mathcal{U}$-code with inner distribution $\bm{b}$ where $\mathcal{T}\cup\mathcal{U}=\{1,2,\dots,n\}$, then by the right side of \eqref{eqn:orthogonality} we have
\begin{equation}\label{eqn:|Y||Z|<=|X|}
	|Y||Z|\le (\overline{\bm{a}Q})\Delta_m^{-1}(\bm{b}Q)^{\mathsf{T}}=|X|\bm{a}\Delta_k^{-1}\bm{b}^{\mathsf{T}}=|X|
\end{equation}
with equality if and only if $(\bm{a}Q)_i(\bm{b}Q)_i=0$ $(1\le i\le n)$.
This ``Anticode Bound'' is a special case of the linear programming method. A similar argument gives an ``Antidesign Bound'' for 
$\mathcal{T}$-designs: if $Y$ is a $\mathcal{T}$-design and $Z$ is a $\mathcal{U}$-design where $\mathcal{T}\cup\mathcal{U}=\{1,2,\dots,n\}$, then $|Y||Z| \ge |X|$.

In some cases, certain $\mathcal{T}$-designs have natural geometric interpretations.
For example, if $(X,\mathcal{R})$ is induced on the top fiber of a short\footnote{A ranked, meet semilattice $(\mathcal{P},\preccurlyeq)$ with top fiber $X$ is \emph{short} if $X\wedge X=\mathcal{P}$.} \emph{regular semilattice} $(\mathcal{P},\preccurlyeq)$ (see \cite{Delsarte1976JCTA}), then $Y$ is a $t$-design\footnote{Here we are using the ordering $\{E_i\}_{i=0}^n$ defined naturally by the semilattice structure.} if and only if the number $|\{x\in Y:u\preccurlyeq x\}|$ (called the \emph{index}) is independent of $u\in\mathcal{P}$ with $\rank(u)=t$ \cite{Delsarte1976JCTA}.
For Examples \ref{exmp:Johnson} and \ref{exmp:Hamming} we have:

\begin{example}\label{exmp:Boolean-lattice} 
Let $\mathcal{P}=\{u\subseteq\Omega:|u|\le n\}$.
Then $(\mathcal{P},\preccurlyeq)$, where the partial order is given by inclusion, forms a short regular semilattice (\emph{truncated Boolean lattice}) with rank function $\rank(u)=|u|$.
In the top fiber $J(v,n)$, a Delsarte $t$-design is just a combinatorial $t$-design.\footnote{A $t$-$(v,n,\lambda)$ \emph{design} is a collection of $n$-subsets (called blocks) of a $v$-set such that every $t$-subset is contained in exactly $\lambda$ blocks.}
\end{example}

\begin{example}\label{exmp:Hamming-lattice} 
Introduce a new symbol ``$\cdot$'' and let $\mathcal{P}$ be the set of words of length $n$ over $\mathcal{Q}\cup\{\cdot\}$.
For $u=(u_1,u_2,\dots,u_n), v=(v_1,v_2,\dots,v_n)\in\mathcal{P}$, we set $u\preccurlyeq v$ if and only if $u_i=\cdot$ or $u_i=v_i$, for all $i$.
Then $(\mathcal{P},\preccurlyeq)$ defines a short regular semilattice (\emph{Hamming lattice}) with rank function $\rank(u)=|\{i:u_i\ne\cdot\}|$.
In the top fiber $H(n,q)$, a $t$-design is an orthogonal array of strength $t$.\footnote{An \emph{orthogonal array} $OA_\lambda(t,n,q)$ is a $\lambda q^t \times n$
matrix over an alphabet $\mathcal{Q}$ of size $q$ in which each set of $t$ columns contains each $t$-tuple over $\mathcal{Q}$ exactly
$\lambda$ times as a row.}
\end{example}

See \cite{Munemasa1986GC,Stanton1986GC} for geometric interpretations of $t$-designs in other classical families of metric \& cometric association schemes.
More ``exotic'' types of codes and designs will be discussed in \S \ref{sec:ubiquity}.
See also \cite{CD1993TAMS,Delsarte2004EJC} for another approach to the regularity of $\mathcal{T}$-designs in $J(v,n)$ and $H(n,q)$ in terms of their $t$-\emph{form spaces}.

The \emph{outer distribution} of $Y$ is the $|X|\times(n+1)$ matrix
\begin{equation}
	B=[A_0\chi,A_1\chi,\dots,A_n\chi].
\end{equation}
We also recall the following \emph{four fundamental parameters} of $Y$:
\begin{gather}
	\delta=\min\{i\ne 0:a_i\ne 0\},  \quad \delta^*=\min\{i\ne 0:(\bm{a}Q)_i\ne 0\}, \\
	s=|\{i\ne 0:a_i\ne 0\}|, \quad s^*=|\{i\ne 0:(\bm{a}Q)_i\ne 0\}|.
\end{gather}
We call $\delta,\delta^*,s,s^*$ the \emph{minimum distance}, \emph{dual distance}, \emph{degree} and \emph{dual degree} of $Y$, respectively.\footnote{We also refer to $\tau=\delta^*-1$ as the (\emph{maximum}) \emph{strength} of $Y$.}
These are related with $|Y|$ as follows:

\begin{theorem}[\cite{Delsarte1973PRRS}]\label{thm:P-inequalities}
Suppose that $(X,\mathcal{R})$ is metric with $P$-polynomial ordering $\{A_i\}_{i=0}^n$.
Then $\delta\le 2s^*+1$ and
\begin{equation}\label{eqn:P-inequalities}
	\sum_{i=0}^{\lfloor(\delta-1)/2\rfloor}k_i\le\frac{|X|}{|Y|}\le\sum_{i=0}^{s^*}k_i.
\end{equation}
If $\delta\ge 2s^*-1$ then $Y$ is completely regular, i.e., the $x^{\text{th}}$ row of $B$ depends only on $\partial(x,Y)=\min\{i:B_{xi}\ne 0\}$.
\end{theorem}

\begin{theorem}[\cite{Delsarte1973PRRS}]\label{thm:Q-inequalities}
Suppose that $(X,\mathcal{R})$ is cometric with $Q$-polynomial ordering $\{E_i\}_{i=0}^n$.
Then $\delta^*\le 2s+1$ and
\begin{equation}\label{eqn:Q-inequalities}
	\sum_{i=0}^{\lfloor(\delta^*-1)/2\rfloor}m_i\le |Y|\le\sum_{i=0}^sm_i.
\end{equation}
If $\delta^*\ge 2s-1$ then $(Y,\mathcal{R}^Y)$ is a cometric subscheme with $s$ classes.
\end{theorem}

The inequality in the left side in \eqref{eqn:P-inequalities} (resp.  \eqref{eqn:Q-inequalities}) is the \emph{sphere-packing bound} (resp. \emph{Fisher-type inequality}), and $Y$ is a \emph{perfect code} (resp. \emph{tight design}) if it satisfies equality.
It follows that $Y$ is a perfect code (resp. tight design) if and only if $\delta=2s^*+1$ (resp. $\delta^*=2s+1$).
We remark that the codes with $\delta\in\{2s^*-1,2s^*\}$ in Theorem \ref{thm:P-inequalities} are precisely the \emph{uniformly packed codes} \cite[p.~348]{BCN1989B}.
Completely regular codes have been actively studied because of their importance in the theory of distance-regular graphs; see \cite[Chapter 11]{BCN1989B} and \cite{Martin2004P}.

Suppose now that $(X,\mathcal{R})$ is metric with $P$-polynomial ordering $\{A_i\}_{i=0}^n$.
Pick any $x\in X$ and set $Z_x=\{y\in X:(x,y)\in\bigcup_{i=0}^eR_i\}$, where $e={\lfloor (\delta-1)/2\rfloor}$ is the \emph{packing radius} of $Y$.
We may obtain the sphere-packing bound via \eqref{eqn:|Y||Z|<=|X|} with $Z=Z_x$.
Thus, if $Y$ is perfect, then since the characteristic vector $\psi_x=\sum_{i=0}^eA_i\hat{x}$ of $Z_x$ satisfies $\psi_x^{\mathsf{T}}E_j\psi_x=(\sum_{i=0}^eP_{ji})^2m_j|X|^{-1}$, we find $|\{j\ne 0:\sum_{i=0}^eP_{ji}=0\}|=e(=s^*)$.
In other words, using the notation at the beginning of \S \ref{sec:metric-cometric-schemes}, all the zeros of the \emph{Lloyd polynomial} $\sum_{i=0}^ev_i$ must be in $\{P_{11},P_{21},\dots,P_{n1}\}$.
We remark that this ``Lloyd Theorem'' has a dual, so that we also obtain a strong nonexistence condition on tight designs in general cometric association schemes in terms of the \emph{Wilson polynomial} $\sum_{i=0}^ev_i^*$.
See \cite{Delsarte1973PRRS} for the details.\footnote{In fact, one may find an analogue of this
Wilson polynomial in schemes which are not necessarily cometric. Let
$\mathcal{E}$, $\mathcal{F} \subseteq \{0,1,\ldots,n\}$. Define  
$ \mathcal{E} \star \mathcal{F}$ to be the set of eigenspace indices $k$ ($0\le k\le n$)
such that $q_{ij}^k \neq 0$ for some $i\in \mathcal{E}$ and some $j\in \mathcal{F}$.
Then, if $Y$ is a $\mathcal{T}$-design and $\mathcal{E}$ satisfies $\mathcal{E} \star \mathcal{E}
\subseteq \mathcal{T} \cup \{0\}$, we obtain the lower bound 
$|Y|\ge \sum_{j \in \mathcal{E}} m_j $. See \cite{Martin2001P} for details and further
conditions on the case when equality holds.}

We may derive a lot more structural information on $Y$ by just looking at the four parameters.
For example, it follows from $BQ=|X|[E_0\chi,E_1\chi,\dots,E_n\chi]$ that $\rank(B)=s^*+1$.
Hence, if $(X,\mathcal{R})$ is metric with $P$-polynomial ordering $\{A_i\}_{i=0}^n$, then the \emph{covering radius} $\rho=\max\{\partial(x,Y):x\in X\}$ of $Y$ must satisfy $\rho\le s^*$.
(The right side of \eqref{eqn:P-inequalities} follows from this.)
We call $Y$ \emph{regular} if $\hat{x}^{\mathsf{T}}A_i\chi$ is independent of $x\in Y$ (and thus equals $a_i$) for all $i$.
It is known that

\begin{theorem}[\cite{Delsarte1973PRRS}]
Suppose that $(X,\mathcal{R})$ is metric with $P$-polynomial ordering $\{A_i\}_{i=0}^n$.
If $\delta\ge s^*$ then $Y$ is regular.
\end{theorem}

\begin{theorem}[\cite{Delsarte1973PRRS}]
Suppose that $(X,\mathcal{R})$ is cometric with $Q$-polynomial ordering $\{E_i\}_{i=0}^n$.
If $\delta^*\ge s$ then $Y$ is regular.
\end{theorem}

We refer the reader to \cite{DL1998IEEE,Camion1998B} for more detailed information and the (major) progress up to 1998.
We remark that Delsarte's theory of codes and designs (in metric/cometric association schemes) based on the linear programming method has been naturally extended to various compact symmetric spaces with rank one,\footnote{In this case, the corresponding orthogonal polynomials are Jacobi polynomials \cite[\S 1.8]{KS1998R}.} such as spheres $S^n=SO(n+1)/SO(n)$; see \cite{BB-ejc} for a survey on this topic.

In 2003, Brouwer et al.~\cite{BGKM2003JCTA} introduced the following parameters for $Y$:
\begin{equation}
	w=\max\{i:a_i\ne 0\}, \quad w^*=\max\{i:(\bm{a}Q)_i\ne 0\}.
\end{equation}
We call $w$, $w^*$ the \emph{width} and \emph{dual width} of $Y$, respectively.
They obtained the following results:

\begin{theorem}[\cite{BGKM2003JCTA}]\label{thm:width}
Suppose that $(X,\mathcal{R})$ is metric with $P$-polynomial ordering $\{A_i\}_{i=0}^n$.
Then $w\ge n-s^*$.
If $w=n-s^*$ then $Y$ is completely regular.
\end{theorem}

\begin{theorem}[\cite{BGKM2003JCTA}]\label{thm:dual-width}
Suppose that $(X,\mathcal{R})$ is cometric with $Q$-polynomial ordering $\{E_i\}_{i=0}^n$.
Then $w^*\ge n-s$.
If $w^*=n-s$ then $(Y,\mathcal{R}^Y)$ is a cometric subscheme with $s$ classes.
\end{theorem}
The above results are in contrast with the bounds on $\delta$ and $\delta^*$ in Theorems \ref{thm:P-inequalities} and \ref{thm:Q-inequalities}.
See \cite{BGKM2003JCTA,HS2007EJC} for many interesting examples attaining the bounds in Theorems \ref{thm:width} and \ref{thm:dual-width}.

Suppose now that $(X,\mathcal{R})$ is metric with $P$-polynomial ordering $\{A_i\}_{i=0}^n$ and cometric with $Q$-polynomial ordering $\{E_i\}_{i=0}^n$.
Since $w\ge s$ and $w^*\ge s^*$ we have
\begin{equation}\label{eqn:w+w*>=n}
	w+w^*\ge n.
\end{equation}
If $(X,\mathcal{R})$ is induced on the top fiber of a short regular semilattice $(\mathcal{P},\preccurlyeq)$, then for every $u\in\mathcal{P}$ the subset $Y_u=\{x\in X:u\preccurlyeq x\}$ satisfies $w=n-\rank(u)$ and $w^*=\rank(u)$.\footnote{The characteristic vectors of the $Y_u$ with $\rank(u)=\ell$ span $\sum_{i=0}^{\ell}V_i$ $(0\le \ell\le n)$; see \cite{Delsarte1976JCTA}.}
It is shown in \cite{BGKM2003JCTA} that any code with $w+w^*=n$ in $J(v,n)$ and $H(n,q)$ is isomorphic to a code of the form $Y_u$.
This result was later extended to their $q$-analogues in \cite{Tanaka2006JCTA}.
It should be remarked that for these examples the $Y_u$ again afford metric \& cometric association schemes which belong to the same family as the original.\footnote{At the algebraic level, this is explained from the results in \cite{Tanakapre}.}
This ``hierarchical structure'' appears to be a subject ripe for further investigation.

The \emph{Erd\H{o}s--Ko--Rado Theorem} \cite{EKR1961QJMO} states that for each integer $t$ such that $v>(t+1)(n-t+1)$ the largest codes satisfying $w\le n-t$ in $J(v,n)$ are the $Y_u$ with $\rank(u)=t$.
In fact, the original proof in \cite{EKR1961QJMO} based on the ``shifting technique'' establishes the conclusion under the stronger hypothesis $v\ge t+(n-t)\binom{n}{t}^3$, and the best possible bound $v>(t+1)(n-t+1)$ was obtained in \cite{Wilson1984C} as an application of Delsarte's linear programming method.
The observation that the largest (or extremal) codes in the Erd\H{o}s--Ko--Rado Theorem are those codes satisfying $w+w^*=n$ led to the ``$q$-versions'' of the theorem in full generality; see \cite{Tanaka2006JCTA}.\footnote{This is a consequence of the previous work \cite{FW1986JCTA,Huang1987DM}, together with the classification of codes with $w+w^*=n$. The construction of \emph{Singleton systems} \cite{Delsarte1978JCTA} (i.e., $t$-designs with index one) in bilinear forms schemes plays an important role in the proof (in view of \eqref{eqn:|Y||Z|<=|X|}); see also \cite[p.~192]{Huang1987DM}.}

\section{Duality}\label{sec:duality}

Suppose that $(X,\mathcal{R})$ is a commutative association scheme and that $X$ is endowed with the structure of an abelian group (written multiplicatively) with identity element $1$.
We call $(X,\mathcal{R})$ a \emph{translation association scheme} if for all $0\le i\le n$ and $z\in X$, $(x,y)\in R_i$ implies $(xz,yz)\in R_i$.
This concept is equivalent to that of a \emph{Schur ring} on an abelian group; see \cite{MP-ejc} for a survey on Schur rings.

Let $X^*$ be the character group of $X$.
To each $\varepsilon\in X^*$ we associate the vector $\hat{\varepsilon}=|X|^{-1/2}\sum_{x\in X}\overline{\varepsilon(x)}\hat{x}$, so that $\langle\hat{x},\hat{\varepsilon}\rangle=|X|^{-1/2}\varepsilon(x)$.\footnote{See footnote \ref{footnote:idempotents}.}
Note that the $\hat{\varepsilon}$ form an orthonormal basis for $\mathbb{C}^X$ by the orthogonality relations for the characters.
Define a partition $\mathcal{X}=\{X_0,X_1,\dots,X_n\}$ of $X$ by $X_i=\{x\in X:(1,x)\in R_i\}$ $(0\le i\le n)$.\footnote{Such a partition is sometimes referred to as a \emph{blueprint}; see e.g., \cite{Bailey2004B}.}
Then $R_i=\{(x,y)\in X\times X:yx^{-1}\in X_i\}$ $(0\le i\le n)$ and we find
\begin{equation}\label{eqn:eigenvectors}
	A_i\hat{\varepsilon}=\left(\sum_{x\in X_i}\overline{\varepsilon(x)}\right)\hat{\varepsilon} \quad (0\le i\le n,\ \varepsilon\in X^*).
\end{equation}
Hence we may also partition $X^*$ as follows: $\mathcal{X}^*=\{X_0^*,X_1^*,\dots,X_n^*\}$, where $X_i^*=\{\varepsilon\in X^*:\hat{\varepsilon}\in V_i\}$ $(0\le i\le n)$.
It follows that
\begin{equation}\label{eqn:eigenvalues}
	P_{ij}=\sum_{x\in X_j}\overline{\varepsilon(x)} \quad (\varepsilon\in X_i^*),  \quad Q_{ij}=\sum_{\varepsilon\in X_j^*}\varepsilon(x)  \quad (x\in X_i)
\end{equation}
for $0\le i,j\le n$.
The left-hand equation of \eqref{eqn:eigenvalues} is immediate from \eqref{eqn:eigenvectors}, and the right-hand 
equation follows by evaluating $|X|(E_j)_{1x}$ in two ways using $E_j=\sum_{\varepsilon\in X_j^*}\hat{\varepsilon}\overline{\hat{\varepsilon}}\,{}^{\mathsf{T}}$.\footnote{Note that $\{\hat{\varepsilon}:\varepsilon\in X_i^*\}$ forms an orthonormal basis for $V_i$.}
Let $\mathcal{R}^*=\{R_0^*,R_1^*,\dots,R_n^*\}$ be the partition of $X^*\times X^*$ defined by $R_i^*=\{(\varepsilon,\eta):\eta\varepsilon^{-1}\in X_i^*\}$ $(0\le i\le n)$, and let $A_i^*\in \mathbb{C}^{X^*\times X^*}$ be the adjacency matrix of $(X^*,R_i^*)$ $(0\le i\le n)$.
If we identify $\mathbb{C}^{X^*\times X^*}$ with $\mathbb{C}^{X\times X}$ via the orthonormal basis $\{\hat{\varepsilon}:\varepsilon\in X^*\}$, then it follows from the orthogonality relations and \eqref{eqn:eigenvalues} that
\begin{equation}\label{eqn:Aj*-in-Ei*}
	A_i^*=\sum_{(\varepsilon,\eta)\in R_i^*}\hat{\varepsilon}\overline{\hat{\eta}}\,{}^{\mathsf{T}}=\sum_{j=0}^nQ_{ji}E_j^*,
\end{equation}
where $E_i^*$ is the diagonal matrix in $\mathbb{C}^{X\times X}$ with $(x,x)$-entry $(E_i^*)_{xx}=(A_i)_{1x}$; so the vector space $\bm{A}^*$ spanned by the $A_i^*$ is an algebra.
Hence $(X^*,\mathcal{R}^*)$ is again a translation association scheme, called the \emph{dual} of $(X,\mathcal{R})$.
By \eqref{eqn:Aj*-in-Ei*}, $(X^*,\mathcal{R}^*)$ has eigenmatrices $P^*=Q$ and $Q^*=P$.
This duality was first formulated in \cite{Tamaschke1962/1963MZ}, but the structure of the Terwilliger algebra (\S \ref{sec:Terwilliger-algebra}) is already visible here.

Let $Y$ be a subgroup of $X$ with characteristic vector $\chi$ and inner distribution $\bm{a}=(a_0,a_1,\dots,a_n)$.
Note that $Y$ is regular and thus $a_i=|Y\cap X_i|$ $(0\le i\le n)$.
Set $Y^{\circ}=\{\varepsilon\in X^*:\varepsilon(y)=1\ \text{for all}\  y\in Y\}$.
Then $Y^{\circ}$ is a subgroup of $X^*$ and
\begin{equation}\label{eqn:MacWilliams-identity}
	(\bm{a}Q)_i=\frac{|X|}{|Y|}\chi^{\mathsf{T}}E_i\chi=\frac{|X|}{|Y|}\sum_{\varepsilon\in X_i^*}|\langle \chi,\hat{\varepsilon}\rangle|^2=|Y|\cdot |Y^{\circ}\cap X_i^*|.
\end{equation}
It follows that $Y^{\circ}$ has inner distribution $|Y|^{-1}\bm{a}Q$.
We remark that $H(n,q)$ is a translation association scheme if we take the alphabet $\mathcal{Q}$ to be an abelian group (cf. Example \ref{exmp:Hamming}).\footnote{The most familiar case is that $\mathcal{Q}$ is a finite field $\mathbb{F}_q$ and $Y$ is a linear code in the usual sense.}
Moreover, in this case the dual of $H(n,q)$ is again the Hamming scheme (with vertex set $X^*=(\mathcal{Q}^*)^n$); in other words, $H(n,q)$ is \emph{self-dual}.
Thus, in view of the generating functions for the Krawtchouk polynomials \cite[\S 1.10]{KS1998R}, \eqref{eqn:MacWilliams-identity} turns out to generalize the well-known MacWilliams identity on the weight distributions (or enumerators) of a linear code and its dual code.
The following theorem is also important (cf. Theorems \ref{thm:Q-inequalities} and \ref{thm:dual-width}):

\begin{theorem}[\cite{Delsarte1973PRRS}]
With the above notation, $(Y,\mathcal{R}^Y)$ is a subscheme if and only if the outer distribution of $Y^{\circ}$ has $s+1$ distinct rows, where $s$ is the degree of $Y$.
\end{theorem}

If $(Y,\mathcal{R}^Y)$ is a subscheme (with $s$ classes), then its dual scheme has vertex set $X^*/Y^{\circ}$ and the relation containing a pair $(\varepsilon Y^{\circ},\eta Y^{\circ})$ is determined by the $(\eta\varepsilon^{-1})^{\rm th}$ row of the outer distribution of $Y^{\circ}$; see \cite{Delsarte1973PRRS}.

Certain dualities of commutative (but not necessarily translation) association schemes also arise in connection with \emph{spin models} and \emph{type II matrices}.
Let $A$ be a nowhere zero matrix in $\mathbb{C}^{X\times X}$ with ``Schur inverse'' $A^{(-)}$, i.e., $A\circ A^{(-)}=J$.
(Henceforth we shall \emph{not} assume a group structure on $X$.)
We call $A$ \emph{type II} if $AA^{(-)\mathsf{T}}=|X|I$.
The \emph{Nomura algebra} of $A$ is the space $\mathcal{N}_A$ of matrices $M$ in $\mathbb{C}^{X\times X}$ such that $A\hat{x}\circ A^{(-)}\hat{y}$ is an eigenvector of $M$ for all $x,y\in X$.
If $A$ is invertible, then $A$ is type II if and only if $J\in\mathcal{N}_A$ (cf. \cite[Lemma 2.1]{CGpre}).
Define a linear map $\Theta_A:\mathcal{N}_A\rightarrow\mathbb{C}^{X\times X}$ by
\begin{equation*}
	M(A\hat{x}\circ A^{(-)}\hat{y})=(\Theta_A(M))_{xy}\cdot (A\hat{x}\circ A^{(-)}\hat{y}) \quad (M\in\mathcal{N}_A,\ x,y\in X).
\end{equation*}
Jaeger et al.~\cite{JMN1998JAC} showed that if $A$ is type II then $\Theta_A(\mathcal{N}_A)=\mathcal{N}_{A^{\mathsf{T}}}$, and\begin{equation}\label{eqn:duality}
	\Theta_{A^{\mathsf{T}}}(\Theta_A(M))=|X|M^{\mathsf{T}}, \quad \Theta_A(MN)=\Theta_A(M)\circ\Theta_A(N)
\end{equation}
for $M,N\in\mathcal{N}_A$.
It follows that if $A$ is type II then both $\mathcal{N}_A$ and $\mathcal{N}_{A^{\mathsf{T}}}$ are the Bose--Mesner algebras of some commutative association schemes, and $\Theta_A$ gives an isomorphism between them which ``swaps'' the ordinary multiplication and $\circ$.

\emph{Spin models} were introduced by V. Jones \cite{Jones1989PJM} as a tool for creating link invariants, and are characterized (up to scalar multiplication) as those type II matrices $A$ satisfying $A\in\mathcal{N}_A$ \cite[Proposition 9]{JMN1998JAC}.
If $A$ is a spin model then in fact we have $\mathcal{N}_A=\mathcal{N}_{A^{\mathsf{T}}}$ and $\Theta_{A^{\mathsf{T}}}=\Theta_A$ \cite[Theorem 11]{JMN1998JAC}; in this case \eqref{eqn:duality} is equivalent to the condition that the corresponding association scheme is \emph{formally self-dual}, i.e., $P=\overline{Q}$ for some orderings of the $A_i$ and the $E_i$; see \cite{BBJ1997JAC}.
In fact, it was shown that $\mathcal{N}_A$ is not just formally self-dual, but is ``strongly hyper-self-dual'' which is defined using the Terwilliger algebra; see \cite{CN2001JAC}.
See e.g., \cite{CW2005JAC,Curtin2007RJ} for more information on the connections to the Terwilliger algebra and Leonard pairs.
Spin models, as well as \emph{four-weight spin models} \cite{BB1995PJM}, have been studied
via the more general but crisp concept of \emph{Jones pairs} \cite{Chan2001D,CGM2003TAMS,CG2004JCTA}.

A formally self-dual association scheme $(X,\mathcal{R})$ is said to satisfy the \emph{modular invariance property} (with respect to $P$) if there is a diagonal matrix $\Delta$ such that $(P\Delta)^3$ is a nonzero scalar matrix.
This gives a necessary condition that $\bm{A}=\mathcal{N}_A$ for a spin model $A$ \cite{BBJ1997JAC,JMN1998JAC}.
The modular invariance property is also quite relevant to \emph{fusion algebras} in conformal field theory \cite{Bannai1993JAC,Gannon2005JAC}.
We remark that fusion algebras are closely related to \emph{table algebras} \cite{Blau-ejc} and to \emph{character algebras} \cite[\S 2.5]{BI1984B} which may in turn be viewed as ``Bose--Mesner algebras at the algebraic level.''

\section{The linear programming bound}\label{sec:LP-bound}

Suppose that $(X,\mathcal{R})$ is a symmetric association scheme.
Let $\mathcal{T}$ be a subset of $\{1,2,\dots,n\}$.
In coding theory, we are often interested in finding a sharp upper bound on the size of a $\mathcal{T}$-code in $X$.
The fact that the inner distribution $\bm{a}$ of a code and its ``MacWilliams transform'' $\bm{a}Q$ are nonnegative leads to the \emph{linear programming} (or \emph{LP}) \emph{bound} developed by Delsarte:

\begin{theorem}[\cite{Delsarte1973PRRS}]
With variable $\bm{a}=(a_0,a_1,\dots,a_n)\in\mathbb{R}^{n+1}$, set
\begin{equation}\label{eqn:LP}
	\ell_{\mathrm{LP}}=\ell_{\mathrm{LP}}(X,\mathcal{T})=\max\ (\bm{a}Q)_0
\end{equation}
subject to (i) $a_0=1$; (ii) $a_i\ge 0$ $(1\le i\le n)$; (iii) $(\bm{a}Q)_{i}\ge 0$ $(1\le i\le n)$; (iv) $a_i=0$ if $i\in\mathcal{T}$.
If $Y\subseteq X$ is a $\mathcal{T}$-code, then $|Y|\le\ell_{\mathrm{LP}}$.
\end{theorem}

The LP bound was shown to be a close variant of Lov\'{a}sz's $\vartheta$-bound \cite{Lovasz1979IEEE} on the Shannon capacity of a graph; see \cite{Schrijver1979IEEE}.
Many computational software packages implement the simplex method to solve linear programming problems, and \eqref{eqn:LP} 
does produce a lot of sharp upper bounds on the size of codes.
However, most analytic results give bounds using the dual linear program:
\begin{theorem}[\cite{Delsarte1973PRRS}]
With variable $\bm{b}=(b_0,b_1,\dots,b_n)\in\mathbb{R}^{n+1}$, set
\begin{equation}\label{eqn:dual-LP}
		\ell_{\mathrm{LP}}'=\ell_{\mathrm{LP}}'(X,\mathcal{T})=\min\ (\bm{b}Q^{\mathsf{T}})_0
\end{equation}
subject to (i) $b_0=1$; (ii) $b_i\ge 0$ $(1\le i\le n)$; (iii) $(\bm{b}Q^{\mathsf{T}})_i\le 0$ if $i\in\{1,2,\dots,n\}-\mathcal{T}$.
Then $\ell_{\mathrm{LP}}=\ell_{\mathrm{LP}}'$.
\end{theorem}

If $\bm{a},\bm{b}$ are feasible solutions to the programs \eqref{eqn:LP} and \eqref{eqn:dual-LP} respectively, then\begin{equation}
	(\bm{a}Q)_0\le\bm{a}Q\bm{b}^{\mathsf{T}}\le (\bm{b}Q^{\mathsf{T}})_0
\end{equation}
with equality if and only if $(\bm{a}Q)_ib_i=a_i(\bm{b}Q^{\mathsf{T}})_i=0$ for $1\le i\le n$.
Note that \eqref{eqn:|Y||Z|<=|X|} also amounts to constructing a feasible solution to \eqref{eqn:dual-LP}; we saw in \S \ref{sec:codes-designs} that the optimality condition was the key to prove and generalize Lloyd's Theorem.
In passing, \eqref{eqn:|Y||Z|<=|X|} can be slightly strengthened as follows: $\ell_{\mathrm{LP}}(X,\mathcal{T})\ell_{\mathrm{LP}}(X,\mathcal{U})\le |X|$, where $\mathcal{T}\cup\mathcal{U}=\{1,2,\dots,n\}$ \cite{Tarnanen1999EJC}.

The LP bound for $\mathcal{T}$-designs is formulated in a totally analogous manner, so 
we omit the details.  This method provides \emph{lower} bounds on the size of designs; 
see \cite{Delsarte1973PRRS}. Due to the divisibility conditions on $|Y|$ inherent in 
the definition of a $\mathcal{T}$-design $Y$ in the familiar cases, these bounds are
not often as sharp as the corresponding bounds for codes, but in many cases, these 
are still the best known non-trivial bounds.\footnote{See, e.g., Table 4.44 in 
\cite{KL2007B} where parameter sets for $t$-$(v,k,\lambda)$ block designs are ruled
out (actually by Haemers, Weug and Delsarte) using linear programming.}

From now on, suppose that $(X,\mathcal{R})$ is the Hamming scheme $H(n,q)$.
The most traditional case here is that $\mathcal{T}$ is of the form 
$\{1,2,\dots, \delta-1\}$ for some $1\le\delta\le n$,
so that we seek an upper bound on $A_q(n,\delta)$, the maximum size of a code in $X$ with minimum distance (at least) $\delta$.
Since $Q_{ij}=K_j(i)$ is a (Krawtchouk) polynomial of degree $j$ in $i$, the dual program can be stated entirely in terms of ``Krawtchouk expansions'' of polynomials:
any polynomial $f=\sum_{j=0}^n b_j K_j$ satisfying (i) $b_0=1$; (ii) $b_j\ge 0$ $(1\le j\le n)$; and (iii) $f(i)\le 0$ $(\delta\le i\le n)$, yields an upper bound $A_q(n,\delta)\le f(0)$.
Hence one may demonstrate feasible solutions for given ranges $\delta$ without necessarily solving to optimality.
Examples of bounds which can be derived in this way are the Plotkin bound\footnote{For $\delta >  (1-1/q)n$, the Plotkin bound uses $f= K_0+b_1 K_1$ where $b_1=1/(q\delta-n(q-1))$ to yield $A_q(n,\delta)\le q\delta/(q\delta-n(q-1))$.} and the bound of McEliece et al.~\cite{MRRW1977IEEE}.
See also \cite{DL1998IEEE} for detailed discussions on the bounds of V. Levenshtein.

In 2001, A. Samorodnitsky \cite{Samorodnitsky2001JCTA} proved that, asymptotically, the optimum solution to Delsarte's LP bound is no better than the average of the upper bound of McEliece et al.~and the Gilbert--Varshamov lower bound:
\begin{equation}
	1 - H_2(\theta) \le R \precsim H_2\left(\frac{1}{2} - \sqrt{\theta(1-\theta)}\right)
\end{equation}
where $H_2(x)=-x\log_2(x)-(1-x)\log_2(1-x)$ is the binary entropy function and $R =\limsup_{n\rightarrow \infty}n^{-1}\log_2A_2(n,\theta n)$ is the asymptotic rate of the largest binary code with ``relative minimum distance'' $\theta$ ($=\delta/n$).
In fact, these two bounds do not coincide for all $0<\theta< 1/2$, so that even if the lower bound obtained from the Gilbert--Varshamov argument is close to the true optimal value of $R$, the linear programming method, alone, will never be sufficient to prove this.

On the other hand, new ideas for obtaining upper bounds on codes are on the horizon.
The theorem of A. Schrijver \cite{Schrijver2005IEEE} applies semidefinite programming to optimize over the positive-semidefinite cone of the Terwilliger algebra of $H(n,2)$; see \S \ref{sec:SDP-bound}.
There have also been attempts to add more constraints to the program defining $\ell_{\mathrm{LP}}$ using geometric arguments; see e.g., \cite[\S IV]{MEL2002IEEE}.
See also \cite{MEL2002IEEE,MEL2007AMC} for another approach which focuses on the ``holes'' of codes in metric association schemes.
It should be remarked that the recent determination of the kissing number\footnote{\label{footnote:kissing-number}The \emph{kissing number} $k(n)$ is the maximum number of unit spheres which can simultaneously touch the unit sphere in $n$-dimensional Euclidean space without pairwise overlapping.} in four dimensions (i.e., $k(4)=24$) \cite{Musin2006DCG,Musin2008AM} is also based on an improvement of the LP bound for spherical codes.
See \cite{BB-ejc}.

\section{The ubiquity of codes and designs}\label{sec:ubiquity}

In this section, we summarize broad applications of Delsarte Theory and the LP bound to various combinatorial objects. For the most part, these examples involve:
an application; an association scheme $(X,\mathcal{R})$; a partially ordered set related to $X$; a partial order on the eigenspaces of $(X,\mathcal{R})$, somehow related to the first partial order.
The short regular semilattices of Examples \ref{exmp:Boolean-lattice} and \ref{exmp:Hamming-lattice} are the paradigmatic examples, but we aim to emphasize more exotic
settings of this same general flavor. In the cometric case, the eigenspaces are naturally ordered linearly; the $E_i$ are
indexed by the elements of the chain (totally ordered set) $\mathcal{C}_n = \{0,1,\ldots,n\}$.
We may also consider products.\footnote{The \emph{product} of two association schemes $(X,\mathcal{R})$, $(Y,\mathcal{S})$ with Bose--Mesner algebras $\bm{A}$, $\bm{B}$ has vertex set $X\times Y$ and Bose--Mesner algebra $\bm{A}\otimes\bm{B}\subseteq\mathbb{C}^{X\times X}\otimes\mathbb{C}^{Y\times Y}=\mathbb{C}^{(X\times Y)\times(X\times Y)}$.}
The eigenspaces of a product of $m$ cometric association schemes are
naturally ordered by a product of chains.
Next, we consider extensions\footnote{The Bose--Mesner algebra of the $m$-fold \emph{extension} \cite{Delsarte1973PRRS} of an association scheme $(X,\mathcal{R})$ is the $m^{\rm th}$ tensor space of that of $(X,\mathcal{R})$. The eigenmatrices of extensions are described in \cite{Tarnanen1987P} using generating 
functions and in \cite{MT2004PAMS} as Aomoto--Gelfand multivariate hypergeometric functions. See also \cite{AM2003JCTA,Mizukawa2004AM}.} of commetric association schemes.
In this case, the partial order to consider on the eigenspaces is the quotient of the $m$-fold
product of chains $\mathcal{C}_n$
under the action of $S_m$.
This poset is obviously isomorphic to a 
downset\footnote{In a partially ordered set $(\mathcal{P},\preccurlyeq)$, a subset $S\subseteq \mathcal{P}$ is a \emph{downset} (or
\emph{lower ideal}) if $x\in S$ and $y\preccurlyeq x$ always imply $y\in S$.} in Young's lattice (all partitions, or Ferrers diagrams, ordered by
inclusion). It turns out that the theory allows us to take further products and extensions of all of these examples to obtain
more.

Since our main goal in this section is to demonstrate the widespread applicability of the theory discussed here, 
we now give  a number of design-theoretic settings where Delsarte's concept of $\mathcal{T}$-design applies. For
each, we describe the combinatorial objects in question,
the association schemes in which they can be found, the relevant partial order on the eigenspaces of these schemes, and -- relative
to this indexing of eigenspaces -- the subset $\mathcal{T}$ for which these objects are Delsarte $\mathcal{T}$-designs. Applications
are discussed in the references.

\begin{example}
For block designs of strength $t$, i.e., $t$-($v,n,\lambda$) designs, the association scheme is the Johnson scheme $J(v,n)$ with
poset $\mathcal{C}_n$ and $\mathcal{T}=\{1,2,\ldots,t\}$. 
\end{example}

\begin{example}
For an orthogonal array of strength $t$,
the ambient association scheme is the  
Hamming scheme $H(n,q)$ with
poset $\mathcal{C}_n$ and $\mathcal{T}=\{1,2,\ldots,t\}$. 
\end{example}

\begin{example}
The incidence graph of a symmetric design is always distance-regular and has two $Q$-polynomial orderings.
In \cite{Martin2001JCMCC}, a number of geometric substructures in finite projective
spaces are shown to be Delsarte $\mathcal{T}$-designs in the corresponding cometric schemes. In all cases,
the poset is $\mathcal{C}_3$ and we have $\mathcal{T}\subseteq \{1,2\}$.
\end{example}

\begin{example}
In \cite{Levenshtein1997DCC}, Levenshtein studies systems of resilient functions with an eye toward cryptographic applications.
The underlying combinatorial objects are  \emph{split orthogonal arrays}, which are Delsarte $\mathcal{T}$-designs
in a product  $H(n_1,q) \otimes H(n_2,q)$.
Here, the eigenspaces are indexed by  $\mathcal{C}_{n_1} \times \mathcal{C}_{n_2}$  and $\mathcal{T}= \mathcal{C}_{t_1} \times \mathcal{C}_{t_2}-\{(0,0)\}$. Levenshtein was the first to derive the LP bound for these objects.
\end{example}

\begin{example}
For some time, statisticians have been using \emph{mixed-level orthogonal arrays} for experimental design, without knowing
whether their constructions were as efficient as they could be. If such an experiment has $n_i$ factors with $q_i$ possible levels
($1\le i\le m$), then, in order to evenly test all $t$-tuples of factors, one seeks a Delsarte $\mathcal{T}$-design in
the product $\bigotimes _{i=1}^m H(n_i,q_i)$. The poset on eigenspaces is
$\times _{i=1}^m \mathcal{C}_{n_i}$ and $\mathcal{T} = \left\{ (j_1,\ldots, j_m) :  0 < \sum_{i=1}^m j_i \le t \right\}$.
The LP bound
for mixed-level orthogonal arrays was derived by Martin \cite{Martin1999DCC} and, independently
and simultaneously, Sloane and Stufken  \cite{SS1996JSPI}.
\end{example}

\begin{example}
One may consider the same set $\mathcal{T}$ for a product of Johnson schemes. An interesting special case
is that of \emph{bipartite block designs} (or \emph{mixed block designs}) \cite{Martin1998JCD}, where points
are colored with two colors, say red and white,  and each block contains $k_1$  red points and $k_2$ white points.
For $i+j\le t$, we require a constant number $\lambda_{i,j}$ of blocks containing any $i$ chosen red points and
any $j$ chosen white ones. If there are $v_1$ red and $v_2$ white points in all, then this is a $\mathcal{T}$-design
in $J(v_1,k_1)\otimes J(v_2,k_2)$ where $\mathcal{T}=\left\{ (i,j): 0< i+j \le t\right\}$.
\end{example}

\begin{example}
A \emph{Room square} of side $n$ is a placement of all the unordered pairs of elements 
from $\Omega=\{1,2,\ldots,n+1\}$ into $\binom{n+1}{2}$ of the cells of an $n\times n$ 
array in such a way that each symbol appears once in each row and once in each column. 
A \emph{Room $d$-cube} of side $n$ is a $d$-dimensional cube of side $n$ in which each 
2-dimensional projection is a Room square of side $n$ \cite{Dinitz2007B}. In the product 
scheme $J(n+1,2)\otimes H(d,n)$, a Room $d$-cube of side $n$ is equivalent to a 
Delsarte $\mathcal{T}$-design \cite{Martin1999DCC}
of (minimal) size $\binom{n+1}{2}$ which is also an $\mathcal{S}$-code
where  $\mathcal{T} = \left\{ (1,0),(2,0),(0,1),(1,1)\right\}$ and 
$\mathcal{S}=\left\{ (1,1), (2,1), \ldots, (1,d-2),(2,d-2) \right\}$.
No example is known in which $d> (n-1)/2$.
\end{example}

\begin{example}
Another recent application of Delsarte's theory of designs, and of the LP bound in
particular, is the discovery of the 
\emph{ordered Hamming scheme}\footnote{Some authors refer to this as ``NRT space,'' after Niederreiter, Rosenbloom and Tsfasman,
whose earlier papers contained some properties of this association scheme without using the association scheme terminology
at all.} \cite{MS1999CJM}.
The most important designs here are the  \emph{ordered orthogonal arrays} (OOAs), which in many cases give rise 
to $(t,m,s)$-\emph{nets} \cite{Niederreiter1987MM,Martin2007B} which
in turn provide quasi-Monte Carlo methods for numerical integration, optimization, and simulation.
For an alphabet $\mathcal{Q}$ of size $q$, form relations $R_0,R_1,\dots,R_\ell$ on $\mathcal{Q}^\ell$ by putting $(x,y)\in R_i$ if their maximal 
common prefix has length $\ell-i$.\footnote{In other words, we consider the $\ell$-fold \emph{wreath product} $H(1,q)\wr\dots\wr H(1,q)$ (see e.g., \cite{Muzychuk2009EJC} for a definition).}
The ordered Hamming scheme $H(s,\ell,q)$ is the $s$-fold extension of the $\ell$-class
symmetric association scheme that results from this construction. Its eigenspaces are ordered by the downset of 
all Ferrers diagrams in Young's lattice that fit inside a rectangle with $s$ rows and $\ell$ columns. OOAs of strength $t$ are
characterized as Delsarte $\mathcal{T}$-designs in $H(s,\ell,q)$ where $\mathcal{T}$ contains all non-empty
Ferrers diagrams with $t$ or fewer cells.
Since $H(s,\ell,q)$ is self-dual, it was natural in \cite{MS1999CJM} also to define ``ordered codes''; these turned out to be 
equivalent to the ``codes for the $m$-metric'' studied in \cite{RT1997PIT}.
See also \cite{MV2007IEEE,Bierbrauer2007DCC}.
\end{example}

We remark that Camion \cite{Camion1998B} also uses extensions of certain commutative association schemes to derive 
MacWilliams identities for various weight enumerators of linear codes, such as complete (or spectral), Lee and 
split weight enumerators, all in a unified manner, together with the results of \S \ref{sec:duality}.

\begin{example}
In the association scheme of the symmetric group $S_n$, a most interesting class of $\mathcal{T}$-designs are the  $\lambda$-\emph{transitive sets} of permutations \cite{MS2006JLMS}. Here the eigenspaces are indexed by partitions of $n$ and the
partial order is reverse dominance order $\unlhd$. A set of permutations is $\lambda$-transitive if it is a Delsarte $\mathcal{T}$-design in this
association scheme, where $\mathcal{T}=\{ \mu : \mu \unlhd \lambda \} - \{ (n)\}$. 
This association scheme also provides a framework for \emph{permutation codes} (or \emph{arrays}), and the corresponding LP bound was studied in detail in \cite{Tarnanen1999EJC}.
These codes are of recent interest because of their application to data transmission over electric power lines; see \cite{CCD2004DCC}.
\end{example}

The concept of ``design systems'' \cite{Martin2001P} is a far-reaching generalization of the poset structures associated with classical metric \& cometric association schemes observed in e.g., \cite{Delsarte1976JCTA,Stanton1986GC}, and establishs a framework which includes all of the above types of Delsarte $\mathcal{T}$-designs as well as the most general bounds for them.\footnote{A \emph{design system} attaches
a poset $(\mathcal{J},\unlhd)$ to the eigenspace indices and embeds the vertex set $X$ in a poset $(\mathcal{P},\preccurlyeq)$ with an
order-preserving surjection $\varphi: (\mathcal{P},\preccurlyeq) \rightarrow(\mathcal{J},\unlhd)$ satisfying three conditions. If $M$
is the incidence matrix of $X$ versus $\mathcal{P}$ (using $\preccurlyeq$ where $X\subseteq \mathcal{P}$), then we require $M$ to have
no repeated columns and, for each $j \in \mathcal{J}$, the submatrix $M_j$ obtained by restricting to columns $x$ with $\varphi(x)=j$
must have constant column sum and column space $W_j$ satisfying $V_j \subseteq W_j \subseteq \oplus_{i\unlhd j} V_i$ where
$V_i$ is the $i^{\mathrm{th}}$ eigenspace of the scheme.}
It should be mentioned that there is another new approach to $\mathcal{T}$-designs based on ``coset geometries'' \cite{Ito2004EJC}.
This approach has the advantage that we can totally forget poset structures (though assuming transitive group actions), so that it may have the possibility to enable more flexible applications.

\section{The Terwilliger algebra}\label{sec:Terwilliger-algebra}

Suppose that $(X,\mathcal{R})$ is a commutative association scheme.
Fix a ``base vertex'' $x\in X$.
For each $i$ ($0\le i\le n$) let $E_i^*=E_i^*(x)$, $A_i^*=A_i^*(x)$ be the diagonal matrices in $\mathbb{C}^{X\times X}$ with $(y,y)$-entries $(E_i^*)_{yy}=(A_i)_{xy}$, $(A_i^*)_{yy}=|X|(E_i)_{xy}$.
Note that $E_i^*E_j^*=\delta_{ij}E_i^*$, $\sum_{i=0}^nE_i^*=I$, and moreover
\begin{equation}
	A_i^*A_j^*=\sum_{k=0}^nq_{ij}^kA_k^*, \quad A_i^*=\sum_{j=0}^nQ_{ji}E_j^*.
\end{equation}
The $E_i^*$ and the $A_i^*$ form two bases for the \emph{dual Bose--Mesner algebra} $\bm{A}^*=\bm{A}^*(x)$ \emph{with respect to} $x$.
The \emph{Terwilliger} (or \emph{subconstituent}) \emph{algebra} $\bm{T}=\bm{T}(x)$ \emph{of} $(X,\mathcal{R})$ \emph{with respect to} $x$ is the subalgebra of $\mathbb{C}^{X\times X}$ generated by $\bm{A}$ and $\bm{A}^*$ \cite{Terwilliger1992JAC,Terwilliger1993JACa,Terwilliger1993JACb}.
The following are relations in $\bm{T}$:
\begin{equation}\label{eqn:relations-in-T}
	E_i^*A_jE_k^*=0 \ \ \text{iff} \ \ p_{ij}^k=0; \quad E_iA_j^*E_k=0 \ \ \text{iff} \ \ q_{ij}^k=0.
\end{equation}
(The latter follows by computing the squared norm of $E_iA_j^*E_k$.)
With the notation of \S \ref{sec:Gelfand-pairs}, we also remark that if $(X,\mathcal{R})=G/K$ where $K$ denotes the stabilizer of $x$ in $G$, then $\bm{T}$ is a subalgebra of the centralizer algebra of $\pi|_K$:
\begin{equation}\label{eqn:T-in-centralizer-algebra}
	\bm{T}\subseteq\{M\in\mathbb{C}^{X\times X}:\pi(g)M=M\pi(g)\ \text{for all}\ g\in K\}.
\end{equation}
Equality in \eqref{eqn:T-in-centralizer-algebra} is known to hold for $H(n,q)=(S_q\wr S_n)/(S_{q-1}\wr S_n)$, for example; see \cite{GST2006JCTA}.

Since $\bm{T}$ is closed under conjugate-transpose, it is semisimple and any two non-isomorphic irreducible $\bm{T}$-modules in $\mathbb{C}^X$ are orthogonal.
Describing the irreducible $\bm{T}$-modules is an active area of research; see e.g. \cite{Terwilliger2005GC,MT2006DM,IT2009EJC} and the references therein.
By \eqref{eqn:relations-in-T} we obtain

\begin{lemma}[\cite{Terwilliger1992JAC}]\label{lem:tridiagonal}
Let $W$ be an irreducible $\bm{T}$-module.
Then the following hold:
\begin{enumerate}
\item If $(X,\mathcal{R})$ is metric with $P$-polynomial ordering $\{A_i\}_{i=0}^n$, then $A_1E_i^*W\subseteq E_{i-1}^*W+E_i^*W+E_{i+1}^*W$ $(0\le i\le n)$, where $E_{-1}^*=E_{n+1}^*=0$.
\item If $(X,\mathcal{R})$ is cometric with $Q$-polynomial ordering $\{E_i\}_{i=0}^n$, then $A_1^*E_iW\subseteq E_{i-1}W+E_iW+E_{i+1}W$ $(0\le i\le n)$, where $E_{-1}=E_{n+1}=0$.
\end{enumerate}
\end{lemma}

An irreducible $\bm{T}$-module $W$ is called \emph{thin} (resp.~\emph{dual thin}) if $\dim E_i^*W\le 1$ (resp.~$\dim E_iW\le 1$) for all $i$.
We remark that $J(v,n)$ and $H(n,q)$ are both \emph{thin}\footnote{This concept is in no way related to the ``thin association schemes''
of Zieschang \cite{Zieschang2005B}.} and \emph{dual thin}, i.e., every irreducible $\bm{T}(x)$-module is thin and dual thin for every $x\in X$.
There are several infinite families of metric \& cometric association schemes which have nonthin irreducible $\bm{T}$-modules 
\cite[Note 6.2]{Terwilliger1993JACb} such as classical forms schemes (e.g., bilinear forms schemes). In such cases, in general, the 
determination of all irreducible $\bm{T}$-modules is yet to be settled, with the notable exception of  
the Doob schemes \cite{Tanabe1997JAC}.
See also \cite{BFK2009EJC}.
The irreducible $\bm{T}$-modules of metric \& cometric association schemes are often studied using the theory of \emph{tridiagonal pairs} \cite{ITT2001P}, these being a generalization of Leonard pairs (\S \ref{sec:metric-cometric-schemes}).
Namely, if $(X,\mathcal{R})$ is both metric and cometric then in view of Lemma \ref{lem:tridiagonal}, $A_1$ and $A_1^*$ act on each irreducible $\bm{T}$-module $W$ as a tridiagonal pair (over $\mathbb{C}$).
We remark that $W$ is thin (and dual thin) if and only if this tridiagonal pair is a Leonard pair.
As of this writing, the classification has been worked out by Ito and Terwilliger \cite{ITpre0807} for the tridiagonal pairs over algebraically closed fields which have the most general ``$q$-Racah'' type.
Their proof involves the representation theory of the quantum affine algebra $U_q(\widehat{\mathfrak{sl}}_2)$.
See also \cite{ITpre0805}.
It was earlier shown \cite{IT2009EJC} that for the forms schemes there are four natural algebra homomorphisms from $U_q(\widehat{\mathfrak{sl}}_2)$ to $\bm{T}$, and that $\bm{T}$ is generated by each of their images together with the center $Z(\bm{T})$.\footnote{These homomorphisms arise from an action of the $q$-\emph{tetrahedron algebra} $\boxtimes_q$ \cite{IT2007CA} on $\bm{T}$, so that their images are actually equal.}
It is also an important and urgent next step to ``pull back'' the above representation-theoretic information to the classification problem of metric \& cometric association schemes.

It would be a reasonable project to apply the progress on the Terwilliger algebra and the tridiagonal pairs to the analysis of codes and designs.
This approach is still in its infancy, but it turns out that we may obtain several interesting results even from the elementary facts about 
$\bm{T}$.
As an example, we discuss the \emph{Assmus--Mattson Theorem} \cite{AM1969JCT}, which gives a criterion as to 
when the supports of the words of a fixed weight $k$ in a \emph{linear} code in $H(n,q)$ form a $t$-design (in $J(n,k)$).\footnote{See \cite{Huber2010B} for detailed discussions on the interaction of error-correcting codes with combinatorial designs.}

Let $W$ be an irreducible $\bm{T}$-module.
We recall the following parameters:\footnote{See footnote \ref{footnote:ordering}. Though we do not (explicitly) use in this paper, but the following are also fundamental in the theory: $d=|\{i:E_i^*W\ne 0\}|-1$ (the \emph{diameter}); $d^*=|\{i:E_iW\ne 0\}|-1$ (the \emph{dual diameter}).}
$r=\min\{i:E_i^*W\ne 0\}$ (the \emph{endpoint});
$r^*=\min\{i:E_iW\ne 0\}$ (the \emph{dual endpoint}).
Set $\bm{1}=\sum_{y\in X}\hat{y}$.
Then $\bm{A}\hat{x}=\bm{A}^*\bm{1}$, which is called the \emph{primary} $\bm{T}$-\emph{module}.
It is thin, dual thin and is the unique irreducible $\bm{T}$-module in $\mathbb{C}^X$ satisfying $r=0$ \emph{or} $r^*=0$.

Suppose $(X,\mathcal{R})$ is cometric with $Q$-polynomial ordering $\{E_i\}_{i=0}^n$.
A vector $\chi\in\mathbb{C}^X$ is a \emph{relative} $t$-\emph{design with respect to} $x$ if $E_i\chi\in\mathbb{C}E_i\hat{x}$ 
for $1\le i\le t$ \cite{Delsarte1977PRR}.
If $(X,\mathcal{R})$ is induced on the top fiber of a short regular semilattice $(\mathcal{P},\preccurlyeq)$, then $\chi$ is a relative $t$-design with respect to $x$ if and only if for each $u\in\mathcal{P}$ with $\rank(u)=t$, $\sum_{y\in X,u\preccurlyeq y}\langle\chi,\hat{y}\rangle$ depends only on $\rank(x\wedge u)$ \cite[Theorem 9.8]{Delsarte1977PRR}.
In \cite{Tanaka2009EJC}, $\bm{T}$ was used to give a new proof of the following analogue of the Assmus--Mattson Theorem:

\begin{theorem}[{\cite[Theorem 8.4]{Delsarte1977PRR}}]\label{thm:Assmus-Mattson-v2}
Suppose that $(X,\mathcal{R})$ is cometric with $Q$-poly\-no\-mi\-al 
ordering $\{E_i\}_{i=0}^n$.
Let $Y\subseteq X$ be a code with characteristic vector $\chi$ and dual distance $\delta^*$.
Set $s_x=|\{i\ne 0:E_i^*\chi\ne 0\}|$.
Then $E_{\ell}^*\chi$ is a relative $(\delta^*-s_x)$-design with respect to $x$ for $0\le\ell\le n$.
\end{theorem}

\begin{proof}
Let $U=(\bm{A}\hat{x})^{\perp}$.
Note that $U$ is the linear span of all irreducible $\bm{T}$-modules in $\mathbb{C}^X$ with dual endpoint $r^*>0$.
Set $S=\{i\ne 0: E_i^*\chi\ne 0\}$. 
Then
\begin{equation*}
	\chi|_U\in\left(\sum_{i=\delta^*}^nE_iU\right)\cap\Biggl(\sum_{j\in S}E_j^*U\Biggr),
\end{equation*}
where $\chi|_U$ denotes the orthogonal projection of $\chi$ to $U$.
Since $A_1^*$ generates $\bm{A}^*$ and takes $s_x(=|S|)$ distinct eigenvalues on $\sum_{j\in S}E_j^*U$, it follows that $\bm{A}^*\chi|_U$ is spanned by $\chi|_U,A_1^*\chi|_U,\dots,(A_1^*)^{s_x-1}\chi|_U$.
Hence by Lemma \ref{lem:tridiagonal} (ii) we find
\begin{equation*}
	\bm{A}^*\chi|_U\subseteq\sum_{i=\delta^*-s_x+1}^nE_iU.
\end{equation*}
This shows $E_i\bm{A}^*\chi\subseteq\mathbb{C}E_i\hat{x}$ for $1\le i\le\delta^*-s_x$, and the proof is complete.
\end{proof}

If the irreducible $\bm{T}$-modules with dual endpoint at most $\delta^*-s_x$ are dual thin in Theorem \ref{thm:Assmus-Mattson-v2}, then the conclusion can in fact be significantly strengthened:
$M\chi$ \emph{is a relative} $(\delta^*-s_x)$-\emph{design with respect to} $x$ \emph{for any} $M\in\bm{T}$.
Note also that by dualizing the above arguments we may get another variant of the Assmus--Mattson Theorem for codes in metric association schemes.
See \cite{Tanaka2009EJC} for the details.
Theorem \ref{thm:Assmus-Mattson-v2} (as well as its dual) does not exactly coincide with the original when applied to $H(n,q)$ with $q>2$.
It is interesting, however, to note that if $(X,\mathcal{R})$ is both metric and cometric then recent results on the \emph{displacement} and \emph{split decompositions} \cite{Terwilliger2005GC} can be successfully used to generalize the \emph{original} version:\footnote{The assumption on the semilattice structure is only for the sake of simplicity; see \cite[Example 5.4]{Tanaka2009EJC}. It is assumed in \cite[Theorem 5.2]{Tanaka2009EJC} that the irreducible $\bm{T}$-modules with endpoint at most $t$ and \emph{displacement} \cite{Terwilliger2005GC} zero are thin, but it follows from the results of \cite{Suzuki2005JAC} that this condition is always satisfied.}

\begin{theorem}[{\cite{Tanaka2009EJC}}]\label{thm:Assmus-Mattson-v3}
Suppose that $(X,\mathcal{R})$ is metric with $P$-polynomial ordering $\{A_i\}_{i=0}^n$ and cometric with $Q$-polynomial ordering $\{E_i\}_{i=0}^n$.
Let $Y\subseteq X$ be a code with characteristic vector $\chi$ and dual distance $\delta^*$.
Set $\delta_x=\min\{i\ne 0:E_i^*\chi\ne 0\}$.
Suppose $t\in\{1,2,\dots,n\}$ is such that for every $1\le r\le t$ we have
\begin{equation*}
	|\{r\le i\le n-r:E_i\chi\ne 0\}|\le\delta_x-r,\quad\text{or}\quad |\{r\le i\le n-r:E_i^*\chi\ne 0\}|\le\delta^*-r.
\end{equation*}
If $(X,\mathcal{R})$ is induced on the top fiber of a short regular semilattice $(\mathcal{P},\preccurlyeq)$, then for each $M\in\bm{T}$, $\sum_{y\in X,u\preccurlyeq y}\langle M\chi,\hat{y}\rangle$ is independent of $u\preccurlyeq x$ with $\rank(u)=t$.
\end{theorem}

If $(X,\mathcal{R})=H(n,q)$ and $x$ is the zero vector $(0,0,\dots,0)$ (where $0\in\mathcal{Q}$), then Theorem \ref{thm:Assmus-Mattson-v3} shows that (the complements of) the supports of the words of fixed weight $k$ in $Y$ form a $t$-design (in $J(n,k)\cong J(n,n-k)$) for every $k$.
In particular, the conclusion of the original Assmus--Mattson Theorem is also true for nonlinear codes as well.\footnote{For example, the $[12,6,6]$ extended ternary Golay code has covering radius three, and it follows from Theorem \ref{thm:Assmus-Mattson-v3} that a coset of weight three support $1$-designs.}

A similar approach was also used in \cite{Tanaka2009EJC} to give a new proof of the \emph{minimum distance bound} \cite{Martin2000DCC} for codes in $H(n,q)$.
We saw in \S \ref{sec:duality} that the MacWilliams identities for the weight enumerator of a linear code can be understood from the duality of $\bm{A}$.
The MacWilliams identities for the \emph{biweight enumerator} \cite{MMS1972IEEE} of a binary linear code can then be proved in terms of $\bm{T}$ for $H(n,2)$; see \cite{Martinpre}.
The \emph{harmonic weight enumerators} of linear codes in $H(n,q)$ and their MacWilliams identities studied in \cite{Bachoc1999DCC,Bachoc1999P} use the harmonic analysis for the group $S_{q-1}\wr S_n$ developed in \cite{Dunkl1976IUMJ,Delsarte1978SIAM}, so that we may view these as closely related to the theme discussed in this section.
See also \cite{CDS1991IEEE,Bachoc1999DCC,Tanabe2001DCC} for other proofs of the Assmus--Mattson Theorem based on harmonic analysis.

We remark that the width and dual width of a code mentioned in \S \ref{sec:codes-designs} is quite compatible with the Terwilliger algebra theory.
For instance, there is a more general approach \cite{Suzuki2005JAC} to the width, based on the Terwilliger algebra \emph{with respect to a code} in metric association schemes.
It is a generalization of the results of \cite{GT2002EJC,Terwilliger2002LAA} on thin irreducible $\bm{T}$-modules with endpoint one, and the width of a code and the \emph{tightness} \cite{JKT2000JAC} of distance-regular graphs\footnote{Tight distance-regular graphs have many interesting combinatorial and geometric properties, one of which is that every local subgraph is strongly regular with certain special nontrivial eigenvalues; see \cite{JKT2000JAC,GT2002EJC}.} are discussed together in the unified context of \emph{tight vectors}; see also \cite{HS2007EJC}.
See \cite[\S 8]{Terwilliger2005GC} for a generalization of Inequality \eqref{eqn:w+w*>=n}.

\section{The semidefinite programming bound}\label{sec:SDP-bound}

Throughout this section, suppose that $(X,\mathcal{R})$ is the binary Hamming scheme $H(n,2)=(S_2\wr S_n)/S_n$, so that $X=\mathcal{Q}^n$ where $\mathcal{Q}=\{0,1\}$.
Let $x=(0,0,\dots,0)$ be the zero vector and write $\bm{T}=\bm{T}(x)$, $E_i^*=E_i^*(x)$ $(0\le i\le n)$.
Recall that $\bm{T}$ coincides with the centralizer algebra of $K=S_n$ acting on $X$.

Let $Y\subseteq X$ be a code.
We consider two subsets $\Pi_1,\Pi_2$ of $G=S_2\wr S_n$ defined by $\Pi_1=\{g\in G:x\in gY\}$, $\Pi_2=\{g\in G:x\not\in gY\}$.
For $i\in\{1,2\}$, let
\begin{equation*}
	M_{\mathrm{SDP}}^i=\frac{1}{|Y|n!}\sum_{g\in\Pi_i}\chi_{gY}(\chi_{gY})^{\mathsf{T}}\in\mathbb{C}^{X\times X}
\end{equation*}
where $\chi_{gY}\in\mathbb{C}^X$ denotes the (column) characteristic vector of $gY$.
Since $\Pi_1,\Pi_2$ are unions of right cosets of $G$ by $K$, it follows that $M_{\mathrm{SDP}}^1,M_{\mathrm{SDP}}^2\in\bm{T}$. 
Moreover, since the $\chi_{gY}(\chi_{gY})^{\mathsf{T}}$ are nonnegative and positive semidefinite, so are $M_{\mathrm{SDP}}^1,M_{\mathrm{SDP}}^2$.
By computing the inner products with the $01$-matrices $E_i^*A_jE_k^*$, we readily obtain
\begin{align*}
	M_{\mathrm{SDP}}^1&=\sum_{i,j,k}\lambda_{ijk}E_i^*A_jE_k^*, \quad M_{\mathrm{SDP}}^2=\sum_{i,j,k}(\lambda_{0jj}-\lambda_{ijk})E_i^*A_jE_k^*,
\end{align*}
where
\begin{equation*}
	\lambda_{ijk}=\frac{|X|}{|Y|}\cdot\frac{|\{(y,y',y'')\in Y^3:(y,y',y'')\ \text{satisfies (*)}\}|}{|\{(y,y',y'')\in X^3:(y,y',y'')\ \text{satisfies (*)}\}|},
\end{equation*}
and condition (*) is defined by
\begin{equation*}
	(y,y')\in R_i,\quad (y',y'')\in R_j,\quad (y'',y)\in R_k. \tag{*}
\end{equation*}
By viewing the $\lambda_{ijk}$ as variables we get the following \emph{semidefinite programming} (or \emph{SDP}) \emph{bound} established by A. Schrijver:
\begin{theorem}[\cite{Schrijver2005IEEE}]
Set
\begin{equation*}
	\ell_{\mathrm{SDP}}=\ell_{\mathrm{SDP}}(n,\delta)=\max\sum_{i=0}^n\binom{n}{i}\lambda_{0ii}
\end{equation*}
subject to (i) $\lambda_{000}=1$; (ii) $0\le\lambda_{ijk}\le\lambda_{0jj}$; (iii) $\lambda_{ijk}=\lambda_{i'j'k'}$ if $(i',j',k')$ is a permutation of $(i,j,k)$; (iv) $\sum_{i,j,k}\lambda_{ijk}E_i^*A_jE_k^*\succcurlyeq 0$; (v) $\sum_{i,j,k}(\lambda_{0jj}-\lambda_{ijk})E_i^*A_jE_k^*\succcurlyeq 0$; (vi) $\lambda_{ijk}=0$ if $\{i,j,k\}\cap\{1,2,\dots,\delta-1\}\ne\emptyset$ (where $\succcurlyeq$ means positive semidefinite).
Then $A_2(n,\delta)\le\ell_{\mathrm{SDP}}$.
\end{theorem}

It is known that semidefinite programs can be approximated in polynomial time within any specified accuracy by interior-point methods; see \cite{Todd2001AN}.
See also \cite[\S 7.2]{Gijswijt2005D} for a discussion on how to ensure that computational solutions do give valid upper bounds on $A_2(n,\delta)$.
While Delsarte's LP bound is a close variant of Lov\'{a}sz's $\vartheta$-bound, Schrijver's SDP bound can be viewed as a variant of an extension of the $\vartheta$-bound based on ``matrix cuts'' \cite{LS1991SIAM}; see also \cite[Chapter 6]{Gijswijt2005D}.
In fact, if we define $\bm{a}=(a_0,a_1,\dots,a_n)$ by $a_i=\lambda_{0ii}\binom{n}{i}$ $(0\le i\le n)$, then the condition that $\bm{a}Q$ is nonnegative is equivalent to the positive semidefiniteness of the matrix $M_{\mathrm{LP}}=\sum_{i=0}^n\lambda_{0ii}A_i$, but since $M_{\mathrm{LP}}=M_{\mathrm{SDP}}^1+M_{\mathrm{SDP}}^2$ this is in turn a consequence of the positive semidefiniteness of $M_{\mathrm{SDP}}^1$ and $M_{\mathrm{SDP}}^2$.
A hierarchy of upper bounds based on semidefinite programming was later proposed in \cite{Laurent2007MP}:
\begin{equation*}
	\ell_+^{(1)}\ge\ell_+^{(2)}\ge\dots\ge\ell_+^{(k)}\ge\dots\ge A_2(n,\delta).
\end{equation*}
It turns out that $\ell_{\mathrm{LP}}=\ell_+^{(1)}\ge\ell_{\mathrm{SDP}}\ge\ell_+^{(2)}$.
Each of the $\ell_+^{(k)}$ can be computed in time polynomial in $n$, but the program defining $\ell_+^{(2)}$ already contains $O(n^7)$ variables.
Two strengthenings of $\ell_{\mathrm{SDP}}$ with the same complexity are also given in \cite{Laurent2007MP}.

The SDP bound was also applied to the problem of finding the stability number of the graph $(X,R_{n/2})$ for even $n$ (known as the \emph{orthogonality graph}) in \cite{KP2007EJC}, where it is shown (among other results) that for $n=16$ the SDP bound gives the exact value $2304$, whereas the LP bound only gives much weaker upper bound $4096$.
This problem arises in connection with quantum information theory \cite{GWT2003IEEE}; see also \cite{GN2008SIAM}.

As $M_{\mathrm{SDP}}^1,M_{\mathrm{SDP}}^2$ are $2^n\times 2^n$ matrices, it is in fact absolutely necessary to simplify the program by explicitly describing the Wedderburn decomposition of the semisimple algebra $\bm{T}$.
The decomposition of $\bm{T}$ (as a centralizer algebra) was worked out in \cite{Dunkl1976IUMJ} in the study of addition theorems for Krawtchouk polynomials, but our discussion below emphasizes the use of $\bm{T}$, based on \cite{Go2002EJC}.

Let $W\subseteq\mathbb{C}^X$ be an irreducible $\bm{T}$-module with endpoint $r$.
Then $W$ has dual endpoint $r$, and there is a basis $\{w_i\}_{i=r}^{n-r}$ for $W$ such that
\begin{equation*}
	w_i\in E_i^*W, \quad A_1w_i=(i-r+1)w_{i+1}+(n-r-i+1)w_{i-1} \quad (r\le i\le n-r)
\end{equation*}
where $w_{r-1}=w_{n-r+1}=0$.
Thus, the isomorphism class of $W$ is determined by $r$.
Moreover, it follows that
\begin{equation*}
	\langle w_i,w_j\rangle=\delta_{ij}\binom{n-2r}{i-r}||w_r||^2 \quad (r\le i,j\le n-r).
\end{equation*}
See \cite{Go2002EJC} for the details.
The actions of the $A_i$ on $W$ may be described from the above information as the $A_i$ are Krawtchouk polynomials in $A_1$, but our argument goes as follows.
For integers $i,k,t$ such that $0\le k\le i\le n$ and $0\le t\le \min\{k,n-i\}$, we recall the following normalization of the dual Hahn polynomials found in \cite{CD1993TAMS}:
\begin{equation*}
	\binom{i}{k}Q_t^{i,k}(\lambda^k(z))=\binom{i}{k-t}\binom{n-i}{t}{}_3F_2\bigg(\!\!\!\begin{array}{c} -t,-z,z-n-1 \\ i-n, -k \end{array}\!\!\!\biggm| 1\biggr),
\end{equation*}
where $\lambda^k(z)=k(n-k)-z(n+1-z)$.
If $i+j+k$ is odd then $E_i^*A_jE_k^*=0$ since $H(n,2)$ is bipartite, so suppose that $i+j+k$ is even.
Then it follows that
\begin{equation*}
	E_i^*A_jE_k^*A_2E_k^*=\beta_{j+2}^{i,k}E_i^*A_{j+2}E_k^*+\alpha_j^{i,k}E_i^*A_jE_k^*+\gamma_{j-2}^{i,k}E_i^*A_{j-2}E_k^*,
\end{equation*}
where $\beta_{j+2}^{i,k}=(t+1)(i+1-k+t)$, $\alpha_j^{i,k}=(k-t)(i-k+t)+t(n-i-t)$ and $\gamma_{j-2}^{i,k}=(k+1-t)(n+1-i-t)$, with $t=(j+k-i)/2$.
Using $2A_2=A_1^2-nI$ we find $E_k^*A_2w_k=\lambda^k(r)w_k$ $(r\le k\le n-r)$.
Combining these facts with the three-term recurrence relation for the $Q_t^{i,k}$ \cite[Theorem 3.1]{CD1993TAMS}, we obtain
\begin{equation*}
	E_i^*A_jw_k=Q_t^{i,k}(\lambda^k(r))E_i^*A_{i-k}w_k=Q_t^{i,k}(\lambda^k(r))\binom{i-r}{i-k}w_i
\end{equation*}
for $r\le k\le i\le n-r$, $0\le j\le n$ such that $i+j+k$ is even, where $t=(j+k-i)/2$.
(The $Q_t^{i,k}$ for $t>\min\{k,n-i\}$ are formally defined by the recurrence relation \cite[Theorem 3.1]{CD1993TAMS}.)
Hence, after orthonormalization of the $w_i$, we get the following algebra isomorphism which preserves the positive-semidefinite cones:
\begin{equation*}
	\varphi:\bm{T}\rightarrow\bigoplus_{r=0}^{\lfloor n/2\rfloor}\mathbb{C}^{(n-2r+1)\times(n-2r+1)}
\end{equation*}	
where the $r^{\rm th}$ block of $\varphi(A_j)$ is the symmetric matrix $(a_{i,k}^{j,r})_{i,k=r}^{n-r}$ given by
\begin{equation*}
	a_{i,k}^{j,r}=a_{k,i}^{j,r}=\begin{cases} Q_{(j+k-i)/2}^{i,k}(\lambda^k(r))\binom{i-r}{i-k}\binom{n-2r}{i-r}^{1/2}\binom{n-2r}{k-r}^{-1/2} & \text{if} \ i+j+k \ \text{even}, \\ 0 & \text{if} \ i+j+k \ \text{odd}, \end{cases}
\end{equation*}
for $r\le k\le i\le n-r$, $0\le j\le n$.
See also \cite{Schrijver2005IEEE,Vallentin2009LAA}.

The SDP bound has also been formulated for binary constant weight codes (i.e., codes in $J(v,n)$) in \cite{Schrijver2005IEEE} and for  nonbinary codes in \cite{GST2006JCTA,Gijswijt2005D}.
The description of the irreducible $\bm{T}$-modules becomes more complicated in this case, but this method turns out to improve the LP bound for many parameters.
It seems to be an important problem to decide whether it is possible or not to establish a suitable SDP bound for $t$-designs in $J(v,n)$ or $H(n,q)$.
The SDP bound for spherical codes was formulated in \cite{BV2008JAMS}; it provides a new proof of $k(3)=12$ and $k(4)=24$.\footnote{See footnote \ref{footnote:kissing-number}.}
See also \cite{BB-ejc}.

\section*{Acknowledgments}

The authors are grateful to Andries Brouwer for helpful comments on graph eigenvalues and to Chris Godsil for a careful reading of the manuscript.
While this work was carried out WJM received support through NSA grant number H98230-07-1-0025.
The work of HT was supported by JSPS Grant-in-Aid for Scientific Research No.~20740003.

%%%%%%%%%%%%%%%%%%%%%%%%%%%%%%%%%%%%%%%%%%%
%%%%%%%%%%%%%%%%%%%%%%%%%%%%%%%%%%%%%%%%%%%
%%%%%%%%%%%%%%%%%%%%%%%%%%%%%%%%%%%%%%%%%%%

\end{document}